\documentclass[10pt]{amsart}

\usepackage[utf8]{inputenc}

\usepackage{colonequals}

\usepackage{url}

\usepackage[pagebackref]{hyperref}

\usepackage{rotating}

\usepackage[small, width=0.9\textwidth]{caption}

\usepackage[all]{xy}
\usepackage{pb-diagram, pb-xy}

\usepackage{asymptote}

\usepackage[dvispnames]{xcolor}

\usepackage{amssymb}

\usepackage{mathrsfs}

\numberwithin{equation}{subsection}

\theoremstyle{plain} 

\newtheorem{introthm}{Theorem} 
\newtheorem{introconj}[introthm]{Conjecture}

\newtheorem{introcor}[introthm]{Corollary}
\newtheorem{thm}[equation]{Theorem} 
\newtheorem{lem}[equation]{Lemma}
\newtheorem{cor}[equation]{Corollary}
\newtheorem{prop}[equation]{Proposition}

\newtheorem*{thm*}{Theorem} 

\theoremstyle{definition}
\newtheorem{defn}[equation]{Definition}
\newtheorem{ex}[equation]{Example}

\newtheorem{rem}[equation]{Remark}

\let\tilde = \widetilde

\newcommand{\bbC}{{\mathbb C}}

\newcommand{\bbG}{{\mathbb G}}

\newcommand{\bbN}{{\mathbb N}}
\newcommand{\bbP}{{\mathbb P}}
\newcommand{\bbQ}{{\mathbb Q}}

\newcommand{\bbZ}{{\mathbb Z}}

\newcommand{\rmM}{{\mathrm M}}
\newcommand{\rmS}{{\mathrm S}}

\newcommand{\cleansupp}{\mathrm{supp}} 

\newcommand{\scrA}{{\mathscr A}}

\newcommand{\scrC}{{\mathscr C}}
\newcommand{\scrD}{{\mathscr D}}

\newcommand{\scrF}{{\mathscr F}}

\newcommand{\scrJ}{{\mathscr J}}

\newcommand{\scrL}{{\mathscr L}}
\newcommand{\scrM}{{\mathscr M}}
\newcommand{\scrN}{{\mathscr N}}
\newcommand{\scrO}{{\mathscr O}}

\newcommand{\scrZ}{{\mathscr Z}}

\newcommand{\Mtilde}{{\,\,\,\,\widetilde{\!\!\!\!\scrM}}}

\newcommand{\an}{\mathit{an}}

\newcommand{\Mhol}{\mathrm{Hol}}
\newcommand{\Vect}{\mathrm{Vect}}

\newcommand{\codim}{\mathrm{codim}}

\newcommand{\Supp}{\mathrm{Supp}}
\newcommand{\Stab}{\mathrm{Stab}}
\newcommand{\Sing}{\mathrm{Sing}}
\newcommand{\CH}{\mathrm{CH}}

\newcommand{\Alt}{\mathrm{Alt}}

\newcommand{\one}{{\mathbf{1}}}

\newcommand{\IC}{{\mathrm{IC}}}

\newcommand{\Ad}{\mathscr{A}\hspace*{-0.15em}d} 

\newcommand{\Gl}{\mathit{Gl}}
\newcommand{\Sp}{\mathit{Sp}}
\newcommand{\SO}{\mathit{SO}}

\newcommand{\Sl}{\mathit{Sl}}
\newcommand{\Pic}{{\mathit{Pic}}}
\newcommand{\Hom}{{\mathrm{Hom}}}
\newcommand{\End}{{\mathrm{End}}}

\newcommand{\Aut}{{\mathit{Aut}}}
\newcommand{\id}{{\mathit{id}}}

\newcommand{\CartierDual}[1]{{{\mathrm{Hom}}(#1, \bbG_m)}}

\newcommand{\scrHom}{{\mathscr{H}\kern -.9pt om}}
\newcommand{\scrExt}{{\mathscr{E}\kern -.9pt xt}}

\DeclareMathOperator{\Rep}{{Rep}}

\DeclareMathOperator{\VB}{VB}

\DeclareMathOperator{\CC}{CC}
\DeclareMathOperator{\cc}{cc}

\DeclareMathOperator{\Coh}{{Coh}}

\newcommand{\DR}{\mathrm{DR}}

\makeatletter
\def\@seccntformat#1{%
  \protect\textup{\protect\@secnumfont
    \ifnum\pdfstrcmp{subsection}{#1}=0 \bfseries\fi
    \csname the#1\endcsname
    \protect\@secnumpunct
  }%
}  
\makeatother

\begin{document}

\title[]{Characteristic cycles and the microlocal geometry of the Gauss map, II}
\author{Thomas Kr\"amer}
\address{Institut f\"ur Mathematik, Humboldt-Universit\"at zu Berlin \newline \hspace*{1em} Unter den Linden 6, 10099 Berlin (Germany)}
\email{thomas.kraemer@math.hu-berlin.de}


\begin{abstract}
We show that for the reductive Tannaka groups of semisimple holonomic $\scrD$-modules on abelian varieties, every Weyl group orbit of weights of their universal cover is realized by a conic Lagrangian cycle on the cotangent bundle. Applications include a weak solution to the Schottky problem in genus five, an obstruction for the existence of summands of subvarieties on abelian varieties and a criterion for the simplicity of the arising Lie algebras.
\end{abstract}

\maketitle
\setcounter{tocdepth}{1} 
\tableofcontents

\setcounter{tocdepth}{2} 

\thispagestyle{empty}

\section*{Introduction}

Holonomic $\scrD$-modules and perverse sheaves play an increasing role in algebraic geometry. On semiabelian varieties they are endowed with a convolution product that leads to a Tannakian description; while on affine tori this has been studied by Gabber, Loeser, Katz, Sabbah, Dettweiler and others from many sides, the case of abelian varieties has emerged only more recently. In this paper we further develop the approach to holonomic $\scrD$-modules on complex abelian varieties in~\cite{KraemerMicrolocal}, relating representation theoretic weights with conic Lagrangian cycles on the cotangent bundle. The main new feature is that we focus not on the Weyl group itself but on its orbits. Consequently working with representation rings of the arising groups allows to include characteristic cycles with arbitrary multiplicities and reveals a natural $\lambda$-ring structure on the group of such cycles, which can often be controlled by computations with Chern-Mather classes. In the second half of the paper we illustrate these new methods with some sample applications: A weak solution to the Schottky problem in genus five, a conjecture of Pareschi and Popa on summands of theta divisors, and a criterion for the simplicity of the arising Lie algebras.
 
\medskip 

Let us now discuss the contents of the paper in some more detail. 
Let $A$ be a complex abelian variety of dimension $g$. As a motivating problem for what follows, suppose that we want to know whether a given subvariety $Z\subset A$ decomposes nontrivially as a sum 
$ 
 Z = X+Y
$
where $X, Y\subset A$ are irreducible subvarieties of positive dimension such that the addition morphism $X\times Y \twoheadrightarrow Z$ is generically finite. Famous examples are theta divisors\footnote{Note that summands of indecomposable theta divisors are nondegenerate, so by~\cite[th.~1]{SchreiederDecomposable} the generic finiteness of the addition morphism $X\times Y \twoheadrightarrow Z$ is automatic in this case.} on the Jacobian of a smooth projective curve or on the intermediate Jacobian of a smooth cubic threefold. A conjecture of Pareschi and Popa~\cite[p.~222]{PP}, in a reformulation by Schreieder~\cite[conj.~19]{SchreiederCurveSummands} who in loc.~cit.~proved it for curve summands, says that the existence of such summands characterizes the locus of Jacobians in the moduli space $\scrA_g$ of principally polarized abelian varieties (ppav's):

\begin{introconj} \label{conj:PP}
Let $(A, \Theta)\in \scrA_g(\bbC)$ be an indecomposable ppav. If $\Theta = X + Y$ decomposes nontrivially as a sum, then
\begin{enumerate} 
\item $(A, \Theta)$ is the Jacobian of a smooth projective curve, or
\smallskip
\item $(A, \Theta)$ is the intermediate Jacobian of a smooth cubic threefold, with $g=5$. 
\end{enumerate}
\end{introconj}

In both cases the only decompositions of the theta divisor are the obvious ones, writing it as a sum of copies of the curve or a difference of two copies of the Fano surface of lines on the threefold. Somewhat surprisingly, without a priori assumptions on the dimension or fundamental classes of the summands as in the work of Casalaina-Martin, Popa and Schreieder~\cite{CMPSGenericVanishing, SchreiederCurveSummands, SchreiederDecomposable}, the only proof we know uses the representation theory of reductive groups and the Tannakian methods developed below, see~\cite{KraemerJacobianThetaSummands}. The idea is to reformulate the existence of a nontrivial decomposition in a sum of subvarieties as a statement in the abelian category $\Mhol(\scrD_A)$ of holonomic $\scrD_A$-modules as follows. Recall that for $\scrM\in \Mhol(\scrD_A)$ its de Rham complex
\[
  \DR(\scrM) \;=\; \Bigl[ \, \scrM^\an \to \Omega^1_A\otimes_{\scrO_A} \scrM^\an \to \Omega^2_A \otimes_{\scrO_A} \scrM^\an \to \cdots \, \Bigr]
\]
is a perverse sheaf, and by the Riemann-Hilbert correspondence any perverse sheaf arises as the de Rham complex of a unique regular holonomic module. For a closed subvariety $Z\subset A$ we denote by $\delta_Z\in \Mhol(\scrD_A)$ the unique regular holonomic module such that
\[
 \DR(\delta_Z) \;\simeq\; \IC_Z[\dim(Z)]
\]
is isomorphic to the perverse intersection complex with support $Z$. If $Z=X+Y$ and the sum map $a: X\times Y \to Z$ is generically finite, then by the decomposition theorem of Beilinson, Bernstein, Deligne and Gabber~\cite{BBD,DCM} we get an embedding as a direct summand
\[
 \delta_Z \; \hookrightarrow \; \delta_X * \delta_Y \;=\; a_\dag(\delta_X\boxtimes \delta_Y)
\]
where $*$ is the convolution product from~\cite{KraemerMicrolocal}. The Tannakian formalism allows to control such summands via representations of certain reductive groups. Let us briefly recall the setup: We work in the quotient $\rmM(A) = \Mhol(\scrD_A)/\rmS(A)$ by the Serre subcategory $\rmS(A)\subset \Mhol(\scrD_A)$ of negligible modules, where a module is said to be {\em negligible} if each of its simple subquotients is stable under translations by a non-zero abelian subvariety of $A$. On the quotient the convolution product descends to a product
\[
 *: \quad \rmM(A) \times \rmM(A) \;\longrightarrow\; \rmM(A)
\] 
with the usual associativity, commutativity and duality properties making~$\rmM(A)$ a rigid abelian tensor category. Here by a~{\em tensor category} we mean a symmetric monoidal $\bbQ$-linear pseudoabelian category. The abelian tensor category from above is Tannakian: For any given $\scrM \in \rmM(A)$, let us denote by~ $\langle \scrM \rangle \subset \rmM(A)$ the smallest rigid abelian tensor subcategory containing it, then we have an equivalence of abelian tensor categories
\[
 \omega: \quad \langle \scrM \rangle
 \; \stackrel{\sim}{\longrightarrow} \; \Rep(G)
\]
with the tensor category of finite-dimensional linear representations of an affine algebraic group $G=G(\scrM)$, see section~\ref{sec:tannakian}. In particular, to any subvariety $Z\subset A$ this attaches 
\begin{itemize}
\item a reductive group $G=G(\delta_Z)$,
\item a representation $\omega(\delta_Z) \in \Rep(G)$.
\end{itemize}
If $Z=X+Y$ decomposes as a sum, then the above embedding $\delta_Z \hookrightarrow \delta_X * \delta_Y$ yields an embedding 
\[
 \omega(\delta_Z) \;\hookrightarrow\; \omega(\delta_X)\otimes \omega(\delta_Y)
 \quad \textnormal{in} \quad 
 \Rep(G(\delta_X\oplus \delta_Y)).
\]
This is a very strong condition on the three representations. One way to exploit it is to look at the adjoint representation. Recall that any object of the quotient category $\rmM(A)$ has a unique representative $\scrM \in \Mhol(\scrD_A)$ which is {\em clean} in the sense that it has no negligible sub- or quotient objects. Let $\Ad_Z \in \Mhol(\scrD_A)$ be the unique clean module such that $\omega(\Ad_Z)$ is the adjoint representation of the group~$G(\delta_Z)$ on its Lie algebra. Then the dimensions $d_X = \dim X$ and $d_Y=\dim Y$ of the summands are bounded as follows: 

\begin{introthm}[=\ref{thm:summand}] \label{introthm:summand}
Let $Z\subset A$ be an irreducible subvariety which is not stable under any translation on the abelian variety. Then for any nontrivial decomposition as a sum of geometrically nondegenerate subvarieties $Z=X+Y$ with $d_X+d_Y = d_Z$ one has \smallskip 
\[
  \min \{ d_X, d_Y\} \;\geq \; 
  \delta \;=\; 
  \tfrac{1}{2} \min 
  \bigl\{ \dim \Supp(\scrM) \mid \scrM \hookrightarrow \Ad_Z,  \; \dim \Supp(\scrM) > 0 \bigr\}. \smallskip 
\]
In particular, there can be no such decomposition if $\delta > \lfloor d_Z/2 \rfloor$.
\end{introthm}

To apply this one only needs to compute a single convolution product: We always have 
$
 \Ad_Z \;\hookrightarrow\; \delta_Z * \delta_{-Z}
$
since
$\mathrm{Lie}(G) \hookrightarrow V\otimes V^\vee$ for any faithful $V\in \Rep(G)$.
For theta divisors on ppav's the resulting bound is sharp on Jacobians of curves or intermediate Jacobians of smooth cubic threefolds, in all other known cases it rules out the existence of summands. Conjecture~\ref{conj:PP} would follow from

\begin{introconj} \label{conj:PP2}
If $(A, \Theta)\in \scrA_g(\bbC)$ is an indecomposable ppav which is not the Jacobian of a curve or the intermediate Jacobian of a cubic threefold, then any summand $\scrM \hookrightarrow \Ad_\Theta$ is either a skyscraper sheaf or has support $A$. 
\end{introconj}

Let us take a look at a few examples. After a translation we may assume $\Theta \subset A$ is symmetric. Then $\omega(\delta_\Theta)$ is a symplectic or orthogonal representation depending on whether $g$ is even or odd~\cite[lemma~2.1]{KrWSchottky}. Very often~$G(\delta_\Theta)$ is the full symplectic or orthogonal group. The reason will become clear later, here we only give a simple example: Let
\begin{align} \nonumber
 S_- &\;=\; \bigl\{ \tbinom{2n}{n} \mid n\notin 2\bbZ \bigr\} \cup \bigl\{2^n \mid n\equiv 1,2 \, \mathrm{mod} \, 4 \bigr\} \cup \bigl\{56 \bigr\} \\ \nonumber
 S_+ &\;=\; \bigl\{ \tbinom{2n}{n} \mid n\in 2\bbZ \bigr\} \cup \bigl\{2^n \mid n\equiv 0,3 \, \mathrm{mod} \, 4 \bigr\} \cup \bigl\{7 \bigr\}
\end{align}
be the dimensions of symplectic minuscule resp.~orthogonal weight multiplicity free  representations of the simple complex Lie algebras other than the standard representations, see tables~\ref{tab:minuscule} and~\ref{tab:wmf} in the appendix. A special case of theorem~\ref{thm:theta-multiplicity-free} then gives

\begin{introthm} \label{thm:theta_odp}
Let $(A, \Theta)\in \scrA_g(\bbC)$ be a ppav whose theta divisor $\Theta = - \Theta \subset A$ is smooth except for finitely many ordinary double points $e_1, \dots, e_k$. 
\smallskip 
\begin{enumerate}
 \item If $g$ is even with $g!-2k\notin S_-$ then $ G(\delta_\Theta) = \Sp_{g!-2k}(\bbC)$. \smallskip
 \item If $g$ is odd with $g!-k\notin S_+$ and no two double points differ by a torsion point, then
\[
 G(\delta_\Theta) \;=\;
 \begin{cases}
  \SO_{g!-k}(\bbC) & \textnormal{\em for} \quad e_1+\cdots + e_k \;=\; 0, \\
  O_{g!-k}(\bbC) & \textnormal{\em for} \quad e_1+\cdots + e_k \;\neq\; 0.
 \end{cases}
\]
\end{enumerate}
\end{introthm}

The above assumptions on the singularities are only made for simplicity. For example, lemma~\ref{lem:torsion-difference} illustrates how in some cases the condition on torsion points can be dropped, and similarly one can often treat nonisolated singularities. One should compare this with the situation for Jacobians of curves or intermediate Jacobians of cubic threefolds where $G(\delta_\Theta)$ is much smaller~\cite{KraemerThreefolds}. The exceptions that we have taken out already appear for nonhyperelliptic Jacobians of genus $g=4$ where $k=2$. This is also the only possibly missing case in the stratification of the moduli space $\scrA_g$ defined by $G(\delta_\Theta)$ for $g=4$; this stratification will be discussed in section~\ref{sec:stratification}, it refines both the Andreotti-Mayer stratification by the dimension of the singular locus of the theta divisor and the one by the degree of the Gauss map in~\cite{CGS}. For conjecture~\ref{conj:PP2} it is good news that in so many cases the group $G(\delta_\Theta)$ is the full symplectic or orthogonal group, since then the adjoint representation is just the symmetric or exterior square of the standard representation. Thus in section~\ref{sec:adjoint} we will verify the conjecture in many cases:

\begin{introthm}[=\ref{thm:theta-not-a-sum}] \label{introthm:theta-not-a-sum}
Let $(A, \Theta)\in \scrA_g(\bbC)$. If $G(\delta_\Theta)$ is a symplectic group and~$\omega(\delta_\Theta)$ is its standard representation, then
\[
 \dim \Supp(\Ad_\Theta) \;\geq\; g-1,
\]
and hence $\Theta$ cannot be written as a sum of positive-dimensional subvarieties.
\end{introthm}

One reason why the support estimate works so easily in this case is that here the adjoint representation is irreducible. In general, it is a very important problem to find conditions on a clean module $\scrM \in \Mhol(\scrD_A)$ which ensure that the Lie algebra of the group $G(\scrM)$ is simple modulo its center. In theorem~\ref{thm:simple} we will give various criteria which in particular apply to $\scrM = \delta_\Theta$ for many ppav's, in fact we do not know any example of an indecomposable ppav for which the adjoint module~$\Ad_\Theta$ is not simple. Like all results in this paper, our criteria for the simplicity of the adjoint module rely on the notion of characteristic cycles. Recall that the sheaf $\scrD_A$ has a natural order filtration whose associated graded is the symmetric algebra of the tangent bundle. For $\scrM \in \Mhol(\scrD_A)$ the associated graded for any good filtration can be viewed as a coherent sheaf on the cotangent bundle $T^*A$. Its support is independent of the chosen good filtration and is called the characteristic cycle $\CC(\scrM)$~\cite[def.~2.2.2]{HTT}. It is a conic Lagrangian cycle on $T^*A$ and hence of the form
\[
 \CC(\scrM) \;=\; \sum_{Z\subseteq A} n_Z \Lambda_Z
\quad \textnormal{with} \quad n_Z \;\in\; \bbN_0,
\]
where $\Lambda_Z \subset T^*A$ denotes the conormal variety to a closed subvariety $Z\subset A$. As our Tannakian construction kills negligible modules, it only sees part of the above cycle: Let us say that a component $\Lambda_Z$ is~{\em negligible} if its base $Z \subset A$ is stable under translations by a non-zero abelian subvariety, and {\em clean} otherwise. A conic Lagrangian cycle is called negligible resp.~clean if all components of its support are so. Thus any conic Lagrangian cycle on splits uniquely as a sum of a negligible and a clean part. We denote\footnote{In definition~\ref{def:cc} this notation will be used in a more general framework where clean cycles need not be conic Lagrangian, but for this introduction the present setup suffices.} by
\[
 \cc(\scrM) \;\;=\; \sum_{\Lambda_Z \; \textnormal{clean}} n_Z \Lambda_Z
\]
the clean part of the characteristic cycle. In section~\ref{sec:indexformula} we attach to each $\Lambda_Z$ a Gauss map generalizing the Gauss map for divisors on abelian varieties. We denote its degree by $\deg(\Lambda_Z) \in \bbN_0$ and extend the degree additively to cycles. Coming back to the question if the adjoint module is simple, we show in section~\ref{sec:product} that if the group $G(\scrM)$ is isogenous to a product of two groups, then after applying an isogeny on the abelian variety one has a similar decomposition of $\cc(\scrM)$. This is excluded in many cases:

\begin{introthm}[$\subset$ \ref{thm:simple}]
Let $(A, \Theta)\in \scrA_g(\bbC)$. Then the Lie algebra of $G(\delta_\Theta)$ is simple modulo its center in each of the following cases: \smallskip 
\begin{enumerate}
	\item If $\deg(\Lambda_Z) > \tfrac{1}{3} \deg(\cc(\delta_Z))$. \smallskip
	\item If the theta divisor has only isolated singularities. \smallskip
	\item If the image of $\cc(\delta_Z)$ under any isogeny of abelian varieties is reduced. 
\end{enumerate}
\end{introthm}

The main technical point for the above is to introduce a $\lambda$-ring $\scrL(A)$ of clean conic Lagrangian cycles on $T^* A$ in such a way that on the Grothendieck ring of representations of the Tannakian group~$G=G(\scrM)$, the assignment that sends a holonomic $\scrD_A$-module to the clean part of its characteristic cycle is a $\lambda$-ring homomorphism
\[
 \cc: \quad K^0(\Rep(G))
 \; \longrightarrow \;
 \scrL(A).
\]
The resulting correspondence between Weyl group orbits of weights for $G=G(\scrM)$ and characteristic cycles is based on the microlocalization from~\cite{KraemerMicrolocal}, but with three important novelties: First, rather than considering the monodromy of Gauss maps we only use that the representation ring of a connected reductive group is the Weyl group invariant part of the character ring; this allows to treat cycles with multiplicities. Next, we emphasize $\lambda$-rings to control tensor constructions such as exterior powers or more general Schur functors, where our methods are inspired by the work of Capell, Maxim, Sch\"urmann, Shaneson and Yokura in~\cite{CMSSY}. Finally, we discuss the problem whether an abstract representation of a reductive group is realized by a $\scrD_A$-module if some tensor construction of it is so. For our applications we want tensor constructions in several variables:~By a {\em tensor construction} in~$r$ variables on a tensor category $\scrC$ we mean any functor $S: \scrC^r \to \scrC$ that can be obtained as a composition of direct sums, tensor products and Schur functors. The same definition works in any $\lambda$-ring, see section~\ref{sec:schur}. We then have:

\begin{introthm}[$\subset$ \ref{thm:inverse-galois}]  \label{introthm:inverse-galois}
Let $\scrM \in \Mhol(\scrD_A)$. Let $\widehat{G}$ be an abstract connected reductive group and
\[
 p: \quad \widehat{G} \;\twoheadrightarrow \; G\;=\; G(\scrM)
\]
an isogeny. Let $S$ be a tensor construction in $r$ variables. If $U_1, \dots, U_r\in \Rep(\widehat{G})$ are representations with
\[
 p^*(\omega(\scrN)) \;\simeq\; S(U_1, \dots, U_r) \smallskip
\]
for some $\scrN \in \langle \scrM \rangle = \Rep(G)$, then there exist clean cycles $\Lambda_1, \dots, \Lambda_r\in \scrL(A)$ with
\[
[e]_* \cc(\scrN) \;=\; S(\Lambda_1, \dots, \Lambda_r)
\quad \textnormal{\em for} \quad
e \;=\; \deg(p).
\]
\end{introthm}

We emphasize that a priori the group $\widehat{G}$ need not be the Tannakian group of any holonomic $\scrD_A$-module and so the above is an inverse Galois problem for these groups. For instance, it would be interesting to study the spin covers of orthogonal groups arising for theta divisors. Another case is the Schottky problem; we here only state the nonhyperelliptic case and refer to section~\ref{sec:schottky} for the hyperelliptic case:~Let us say that a ppav $(A, \Theta) \in \scrA_g(\bbC)$ is a {\em nonhyperelliptic fake Jacobian} if the pair $(G(\delta_\Theta), \omega(\delta_\Theta))$ looks like the one for the theta divisor on the Jacobian of a nonhyperelliptic curve. We want to understand whether any fake Jacobian is indeed the Jacobian of such a curve $C \subset A$. For Jacobians $\delta_\Theta = \Alt^{*(g-1)}(\delta_C)$ is an exterior convolution power. Reconstructing the curve from the theta divisor amounts to reversing the effect of this tensor construction. So  theorem~\ref{introthm:inverse-galois} tells us where to look for the curve also in the case of fake Jacobians:

\begin{introthm}[$\subset$ \ref{thm:fake-jacobian}] \label{thm:fake-curve}
If $(A, \Theta) \in \scrA_g(\bbC)$ is a nonhyperelliptic fake Jacobian, then we have
\[ [n]_*\cc(\delta_\Theta) \;=\; \Alt^{g-1}(\Lambda)
\quad \textnormal{\em for some effective cycle $\Lambda \in \scrL(A)$ and some $n\in \bbN$}.
\]
\end{introthm}

When applied to Jacobians this gives a constructive Torelli theorem. For fake Jacobians the cycle $\Lambda$ might however be supported not over a curve but over a higher-dimensional subvariety. Trying to exclude this, we find in section~\ref{sec:weak-schottky} the following solution to the weak Schottky problem in genus five:

\begin{introcor}[=\ref{cor:schottky}]
For $g=5$ the locus of fake Jacobians is contained in the Andreotti-Mayer locus 
\[ \scrN_1 \;=\; \{ (A, \Theta) \in \scrA_5 \mid \dim \Sing(\Theta) \geq 1 \} \;\subset\; \scrA_5, \]
hence it contains the locus of true Jacobians of curves as an irreducible component.
\end{introcor}

The proof is a computation with Chern-Mather classes via Vogel's intersection algorithm~\cite{Vogel,FlennerJoins,FlennerCarrollVogel}, inspired by~\cite{CGS}. Since these techniques are also useful in other contexts, we gather the relevant properties of Chern-Mather classes in section~\ref{sec:chernmather}; in particular, we show that the total Chern-Mather class gives a ring homomorphism from certain subrings of $\scrL(A)$ to truncations of the Chow ring $\CH_\bullet(A)$ endowed with the Pontryagin product. It seems a very interesting problem to understand excess intersections related to positive-dimensional fibers of the Gauss map.

\section{The ring of clean cycles} \label{sec:weights}

We now set up our basic framework for the discussion of characteristic cycles by introducing a~$\lambda$-ring of clean cycles on the cotangent bundle. Let $A$ be a complex abelian variety of dimension $g\geq 1$ and $V=H^0(A, \Omega_A^1)$.

\subsection{Motivation: Kashiwara's formula} \label{sec:indexformula}

The cotangent bundle $T^*A = A\times V$ is canonically trivialized using translations on the abelian variety. We then define the {\em Gauss map} of a subvariety $\Lambda \subset T^* A = A\times V$ of pure dimension $g=\dim(A)$ to be the projection
\[
 \gamma_\Lambda: \quad \Lambda \;\hookrightarrow\; T^* A = A\times V \;\twoheadrightarrow\; V
\]
to the second factor. We will mostly be interested in the case where $\Lambda \subset T^*A$ is conic, i.e.~invariant under rescaling in the fibers of the cotangent bundle. We then write
$ \gamma_{\,\bbP \Lambda}: \bbP \Lambda \rightarrow  \bbP V $
for the projectivized Gauss map. If~$\Lambda$ is taken to be the conormal variety to a divisor $Z\subset A$, this recovers the classical notion of Gauss maps on abelian varieties. Coming back to the general case, for dimension reasons the map $\gamma_\Lambda$ is either generically finite or non-dominant. We say that $\Lambda$ is~{\em clean} in the first case and {\em negligible} in the second case. In what follows we denote the generic degree of the Gauss map by
$\deg(\Lambda) = \deg(\gamma_\Lambda) \in \bbN_0$,
so~$\Lambda$ is negligible iff $\deg(\Lambda)=0$. We extend the degree additively to cycles of pure dimension~$g$ on the cotangent bundle $T^*A$. Then for $\scrM \in \Mhol(\scrD_A)$ we can read off the Euler characteristic 
\[\chi(A, \DR(\scrM)) \;=\; \sum_{i\in \bbZ} (-1)^i \dim_\bbC H^i(A, \DR(\scrM))\]
from the characteristic cycle:

\begin{thm} \label{thm:kashiwara}
We have
$\chi(A, \DR(\scrM)) = \deg(\CC(\scrM))\geq 0$.
\end{thm}

{\em Proof.} This is a special case of Kashiwara's index formula and was shown more generally for semiabelian varieties by Franecki and Kapranov~\cite{FKGauss}. \qed

\medskip 


The above suggests that the most interesting part of characteristic cycles are their clean components, so we should take a closer look at the group of clean cycles on $T^*A$. In what follows we will endow it with a natural $\lambda$-ring structure tailored to fit the convolution product on holonomic $\scrD_A$-modules in section~\ref{sec:tannakian-groups} below.

\subsection{A categorification} 

Motivated by the microlocalization of $\scrD_A$-modules, we categorify the notion of clean cycles by considering germs of vector bundles on their support as in~\cite[sect.~4]{KraemerMicrolocal}:
\begin{defn} \label{defn:vb}
For $u\in V$ let~$\VB(A, u)$ be the $\bbC$-linear exact pseudoabelian category of germs 
\[ \alpha=(\scrF_\alpha, U_\alpha), \]
where $\scrF_\alpha \in \Coh(\scrO_{A\times U_\alpha}^\mathit{an})$ is an analytic coherent sheaf and $u\in U_\alpha \subset V$ a Zariski open such that 
\begin{itemize}
\item the support $\Supp(\scrF_\alpha) \subset A\times U_\alpha$ is algebraic, \smallskip 
\item the projection $\gamma_\alpha: \Supp(\scrF_\alpha) \to U_\alpha$ is finite and flat, \smallskip 
\item the direct image $\gamma_{\alpha *}(\scrF_\alpha)$ is a coherent analytic vector bundle.
\end{itemize} 
Morphisms of germs are defined by 
\[
 \Hom_{\VB(A, u)}(\alpha, \beta) \;=\; \mathrm{colim} \; \Hom_{\scrO_{A\times U}^\mathit{an}} (\scrF_\alpha|_{A\times U}, \scrF_\beta|_{A\times U})
\]
where the colimit runs over all Zariski open $U\subset U_\alpha \cap U_\beta$ containing~$u$, and short exact sequences are taken to be those that are represented by short exact sequences in $\Coh(\scrO_{A\times U}^\mathit{an})$ for $U\ni u$ small enough.
\end{defn}

\noindent
Each germ in the above category gives rise to a clean cycle as follows:

\begin{defn} \label{def:support_closure} 
For $\varnothing \neq U\subseteq V$ Zariski open, a {\em clean} cycle on $A\times U$ is a finite formal sum
\[
 Z \;=\; \sum_\Lambda \, m_\Lambda \cdot \Lambda
 \quad \textnormal{with coefficients} \quad m_\Lambda \;\in\;\bbZ,
\]
where $\Lambda \subset A\times U$ ranges over all irreducible Zariski closed subvarieties for which the projection to the second factor
$
 \gamma_\Lambda: \Lambda \longrightarrow U
$
is a generically finite dominant map. Thus a cycle is clean iff the Zariski closure of each irreducible component of its support is a clean subvariety of the cotangent bundle. Put
\begin{align*}
  \scrZ(A) &\;=\; \{\,\textnormal{clean cycles on $T^*A$}\,\}, \\
  \scrZ(A, U) &\;=\; \{ \,\textnormal{clean cycles on $A\times U$}\,\}, 
\end{align*}
then the pullback of cycles via the embedding $j: A\times U \hookrightarrow A\times V$ gives a group isomorphism
\[
 \scrZ(A) \;\stackrel{\sim}{\longrightarrow}\; \scrZ(A, U), 
 \quad 
 \Lambda \;\mapsto\; \Lambda|_U \;=\; j^{-1}(\Lambda) \medskip
\]
whose inverse is the map that sends a clean cycle $Z=\sum_\Lambda m_\Lambda \cdot \Lambda$ on $A\times U$ to its Zariski closure
\[ 
 \overline{Z} \;=\; \sum_{\Lambda} m_\Lambda \cdot \overline{\Lambda}.
\]
on $T^* A$. This Zariski closure does not change if we replace $U\subset V$ by a smaller Zariski open dense subset. Sending the germ of a sheaf to the Zariski closure of its support we get a homomorphism
\[
\cleansupp: \quad K^0(\VB(A, u)) \;\longrightarrow\; \scrZ(A), \quad \alpha \;=\; (\scrF_\alpha, U_\alpha) \;\mapsto\; \overline{\Supp(\scrF_\alpha)}
\]
on the additive Grothendieck group modulo split exact sequences. 
\end{defn}

\subsection{Lambda rings} \label{sec:grothendieck}
The category $\VB(A, u)$ is naturally a tensor category for the product
\[(\scrF_\alpha, U_\alpha) * (\scrF_\beta, U_\beta) \;=\; (\varpi_*\rho^{-1}(\scrF_\alpha \boxtimes \scrF_\beta), U=U_\alpha \cap U_\beta), \]
where $\rho$ and $\varpi$ are the maps in the correspondence of cotangent bundles induced by the addition on the abelian variety:
\[
\qquad\quad
\xymatrix@M=0.8em@R=-0.5em{
	T^* A^2 \; \supseteq \; A^2 \times U_\alpha \times U_\beta & A^2 \times U \ar[l]_-\rho \ar[r]^-\varpi & A\times U \; \subseteq \; T^* A \\
	\qquad \qquad \qquad ((x,y),u,u) & ((x,y), u) \ar@{|->}[r] \ar@{|->}[l] & (x+y, u) \qquad 
}
\]
This definition is adapted to the convolution product on holonomic $\scrD_A$-modules in section~\ref{sec:tannakian-groups}.
To see what it means for characteristic cycles, recall that for any exact tensor category its Grothendieck group modulo split exact sequences is a ring via the tensor product and a $\lambda$-ring\footnote{We use the term {\em $\lambda$-ring} for what is called a {\em special $\lambda$-ring} in loc.~cit.} via exterior powers~\cite[lemma~4.1]{Heinloth}. In proposition~\ref{prop:lambda} we will upgrade $\cleansupp$ to a $\lambda$-ring homomorphism. Any $\lambda$-ring~$R$ has a natural operation by the ring of symmetric functions in countably many formal variables $x_1, x_2, \dots$~\cite[\S I.3]{AtiyahTall}: If $e_i=e_i(x_1,x_2, \dots)$ are the elementary symmetric functions, the map
\[
 \mu_R: \quad \bbZ[e_1, e_2, \dots] \; \longrightarrow\; \mathrm{Maps}(R, R), \quad e_n \; \mapsto \; \lambda^n
\]
is a ring homomorphism when the target is endowed with pointwise addition and multiplication, and for any homomorphism $f: R\to S$ of $\lambda$-rings and any symmetric function $\sigma$ the diagram
\[
\xymatrix@M=0.5em@C=3em@R=2.5em{
 R \ar[r]^-{\mu_R(\sigma)} \ar[d]_-f & R \ar[d]^-f \\
 S \ar[r]^-{\mu_S(\sigma)} & S
}
\]
commutes. For instance the Adams operations on a $\lambda$-ring $R$ can be defined as the maps
\[ 
\Psi^n \;=\; \mu_R(p_n): \quad R\;\longrightarrow\; R
\quad \textnormal{for the power sums} \quad 
p_n \;=\; x_1^n + x_2^n + \cdots
\] 
As the elementary symmetric functions can be expanded as polynomials in power sums with rational coefficients, the Adams operations determine the $\lambda$-operations on any $\lambda$-ring without $\bbZ$-torsion (see also section~\ref{sec:schur}).

\begin{ex} \label{ex:character-ring}
The {\em representation ring} of a complex algebraic group~$G$ is by definition the Grothendieck ring $R(G) = K^0(\Rep(G))$ of the category of its finite dimensional complex algebraic representations. By restricting representations to a multiplicative subgroup $T\subset G$ we get a homomorphism $R(G) \to R(T)=\bbZ[X]$ to the group ring of the character group $X=\CartierDual{T}$. For maximal tori in  connected reductive groups this restriction homomorphism is an isomorphism onto the invariants
$
 R(G) \stackrel{\sim}{\longrightarrow} \bbZ[X]^W
$
of the Weyl group $W=N_G(T)/Z_G(T)$. Since any representation of a torus is a sum of characters and on characters the Adams operation $\Psi^n$ acts by the $n$-th power map $[n]: X\to X, \chi\mapsto \chi^n$~\cite[II.7.4]{BtD}, we can then describe the Adams operations on the representation ring via the following commutative diagrams:
\[
\xymatrix@M=0.5em{
 R(G) \ar@{^{(}->}[r] \ar[d]_-{\Psi^n} & \bbZ[X] \ar[d]^-{[n]_*} \\
 R(G) \ar@{^{(}->}[r] & \bbZ[X]
}
\]
%
\end{ex} 

\begin{ex} \label{ex:cleanconvolution}
By definition of the tensor structure on $\VB(A, u)$, the support map
$
 \cleansupp:  K^0(\VB(A, u)) \rightarrow \scrZ(A)
$ 
in definition~\ref{def:support_closure} is a ring homomorphism for the product 
\[
 \circ: \quad \scrZ(A)\times \scrZ(A) \;\longrightarrow\; \scrZ(A), \quad 
 \Lambda_1 \circ \Lambda_2 \;=\; 
 \overline{
 \varpi_* \rho^{-1} (\Lambda_1 \times \Lambda_2)|_U.
 }
\]
Here $U\subseteq V$ is any Zariski open dense subset over which the support of both~$\Lambda_i$ is finite and flat, and the closure and the pushforward are again taken in the sense of cycles. There is also a natural structure of a $\lambda$-ring on $\scrZ(A)$. To describe it we need some notations: For $\Lambda \in \scrZ(A)$, let $\langle \Lambda \rangle \subset \scrZ(A)$ be the smallest subring which contains~$\Lambda$ and is stable under taking irreducible components of its members; subrings of this form will be called {\em finitely generated}. For~$u\in V(\bbC)$ consider the subgroup
\[
 \Gamma(\Lambda, u) \;=\; \langle a\in A(\bbC) \mid (a, u)\in \Supp (\Lambda) \rangle 
 \;\subset\; A(\bbC).
\]
Since $\Lambda$ is a clean cycle, this group is finitely generated for general $u\in V(\bbC)$, and its isomorphism type jumps at most for $u$ in a countable union of proper Zariski closed subsets of $V(\bbC)$. For $n\in \bbZ$ we write $[n]: A\times V\to A\times V, (a,v)\mapsto (na, v)$ for the multiplication by $n$ on the abelian variety and denote the induced group homomorphism on clean cycles by
$ [n]_*: \scrZ(A) \to \scrZ(A)$.
\end{ex}

\begin{prop}  
\label{prop:lambda} 
With notations as above we have: \smallskip

\begin{enumerate} 
\item There exists a unique structure of a $\lambda$-ring on $\scrZ(A)$ such that for all $u\in V$ the map
\[
 \cleansupp: \quad K^0(\VB(A, u)) \; \longrightarrow \; \scrZ(A)
\]
is a $\lambda$-ring homomorphism. This $\lambda$-ring structure is given by
$\Psi^n = [n]_*$. \medskip
\item For any given $\Lambda$ the subring 
$\langle \Lambda \rangle \subset \scrZ(A)$
is stable under the $\lambda$-operations, and by taking the fiber over a very general closed point $v\in V(\bbC)$ we get an embedding
\[ 
 (-)_v: \quad \langle \Lambda \rangle \;\hookrightarrow\; \bbZ[\Gamma]
\]
as a $\lambda$-subring
in the group ring of the abelian group $\Gamma = \Gamma(\Lambda, v)$.
\end{enumerate}
\end{prop}

{\em Proof.} If $\eta\in V$ denotes the generic point, then for each closed point $u\in V(\bbC)$ we have an exact inclusion functor $\iota: \VB(A, u)\hookrightarrow \VB(A, \eta)$ which preserves the tensor structure and is compatible with supports. This provides a factorization
\[
\xymatrix@M=0.5em@C=1em{
 K^0(\VB(A, u)) \ar[rr]^-\cleansupp \ar[dr]_{\iota_*\quad} && \scrZ(A) \\
 & K^0(\VB(A, \eta)) \ar[ur]_-{\quad\cleansupp} &
}
\]
where $\iota_*$ is a $\lambda$-ring homomorphism, so for the existence and uniqueness statement it suffices to take $u=\eta$. Note that we can verify these statements entirely in terms of the Adams operations because $\scrZ(A)$, being a ring of cycles with integer coefficients, has no~$\bbZ$-torsion. Moreover it will suffice to consider finitely generated subrings: To check whether a given set of ring endomorphisms $\Psi^n: \scrZ(A)\to \scrZ(A)$ are Adams operations of some $\lambda$-ring structure, we must check a set of universal identities that each only involves finitely many elements at a time. 

\medskip

So let us fix $\Lambda \in \scrL(A)$. The subring $\langle \Lambda \rangle \subset \scrZ(A)$ is countable, so the reduced supports of all its members are finite and \'etale over a very general $v\in V(\bbC)$. In particular, distinct irreducible components of members of~$\langle \Lambda \rangle$ cannot meet over such a point~$v$ since then the union of their reduced supports would not be \'etale over that point. Hence taking the fiber of the occuring cycles over~$v$ defines an injective map $(-)_v: \langle \Lambda \rangle \hookrightarrow \bbZ[\Gamma]$ where $\Gamma = \Gamma(\Lambda, v)$, and the definitions easily imply that this is a ring homomorphism. To see that its image is a $\lambda$-subring, let us denote by $\Vect_\Gamma(\bbC)$ the category of finite-dimensional  $\Gamma$-graded complex vector spaces. The full tensor subcategory
\[
 \VB(A, v)_\Lambda \;=\; \{ \, \alpha \in \VB(A, v) \mid \cleansupp(\alpha) \in \langle \Lambda \rangle \, \} \;\subset\; \VB(A, v)
\]
fits in the commutative diagram
\[
\xymatrix@M=0.5em@C=2em@R=2em{
K^0(\VB(A, v)_\Lambda) \ar[r] \ar@{->>}[d]_-{\cleansupp\;} & K^0(\Vect_\Gamma(\bbC)) \ar@{=}[d]\\
\langle \Lambda \rangle \ar@{^{(}->}[r]^-{(-)_v} & \bbZ[\Gamma]
}
\]
where the top row is induced by the functor that sends $\alpha = (\scrF_\alpha, U_\alpha) \in \VB(A, v)_\Lambda$ to the fiber $(\gamma_{\alpha*}(\scrF_\alpha))(v)$ of the corresponding locally free sheaf, with the grading coming from the inclusion
$
 \{ a\in A(\bbC) \mid (a,v)\in \Supp(\scrF_\alpha)\} \subseteq \Gamma = \Gamma(\Lambda, v)
$. In particular, the top row of the diagram is a homomorphism of $\lambda$-rings and hence it follows that 
\[
 \cleansupp: \quad K^0(\VB(A, v)_\Lambda) \;\longrightarrow\; \langle \Lambda \rangle
\]
is a homomorphism of $\lambda$-rings when the target is endowed with the $\lambda$-ring structure induced by the one of the group ring $\bbZ[\Gamma]=K^0(\Vect_\Gamma(\bbC))$. For the latter the Adams operations $\Psi^n$ are induced by the multiplication by $n$ map on the character group $\Gamma$ as we recalled in example~\ref{ex:character-ring}. It follows that on the subring $\langle \Lambda \rangle \subset \bbZ[\Gamma]$ we have
$\Psi^n = [n]_*$
as claimed. \qed

\subsection{Schur functors} \label{sec:schur}

In any exact tensor category~$\scrC$ we have exterior powers or more generally Schur functors~\cite[sect.~1.4]{DelCT}. By a {\em tensor construction} in $r$ variables we mean any functor 
\[
 S: \quad \scrC^r \;=\; \scrC \times \cdots \times \scrC \; \longrightarrow \; \scrC
\]
that can be obtained as a composition of tensor products, direct sums and Schur functors. 
Any tensor construction descends to an operation on $K^0(\scrC)$; in fact we can define tensor constructions more generally as natural operations on the category of~$\lambda$-rings, with Schur functors replaced by the natural operation of the Schur polynomials \cite[p.~43 and sect.~III.3]{KnutsonLambda}. Thus any homomorphism of $\lambda$-rings is compatible with tensor constructions. Once we know the Adams operations, the effect of Schur functors is easy to describe: For any partition $\beta=(\beta_1, \dots, \beta_\ell)$ we have the power sum polynomials 
\[
 p_\beta \;=\; \prod_{i=1}^\ell \Bigl( \sum_j x_j^{\beta_i}\Bigr).
\]
These power sums are a basis for the ring of symmetric functions with rational coefficients, so for any partition $\alpha = (\alpha_1, \alpha_2, \dots)$ the Schur polynomial $s_\alpha$ has an expansion 
\[
 s_\alpha \;=\; \sum_{\beta} m_{\alpha \beta} p_\beta 
 \quad \textnormal{with unique coefficients} \quad m_{\alpha \beta} \; \in \; \bbQ
\]
where $\beta = (\beta_1, \dots, \beta_\ell)$ runs over all partitions of degre $\deg(\beta)=\deg(\alpha)$. It follows that the operation of the Schur polynomial $s_\alpha$ is given in terms of the Adams operations by
\[
 s_\alpha(x) = \sum_{\beta} m_{\alpha \beta} \cdot \Psi^{\beta_1}(x) \cdot \Psi^{\beta_2}(x) \cdots \Psi^{\beta_\ell}(x).
\]
We thus obtain

\begin{lem} \label{lem:schur}
For any $\Lambda \in \scrZ(A)$ and any partition $\alpha$ we have in $\scrZ(A)\otimes_\bbZ \bbQ$ the formula
\[
 s_\alpha(\Lambda) \;=\; \sum_{\beta} m_{\alpha \beta} \cdot \Lambda_{[\beta]}
 \quad \textnormal{where} \quad
  \Lambda_{[\beta]} \;=\; [\beta_1]_*(\Lambda) \circ \cdots \circ [\beta_\ell]_*(\Lambda)
\]
\end{lem}

{\em Proof.} Immediate since by proposition~\ref{prop:lambda} the Adams operations on the ring of clean cycles are given by $\Psi^n = [n]_*$.
\qed

\medskip 

\begin{ex} \label{ex:elementary-symmetric}
The elementary symmetric polynomials $e_n = s_{1,1,\dots, 1}$ can be written in terms of the power sum polynomials $p_\beta$ from above by expanding the generating series 
\begin{align*}
 \sum_{n=0}^\infty e_n X^n &\;=\; \exp\Bigl( \, \sum_{\nu=1}^\infty \tfrac{(-1)^{\nu + 1}}{\nu} \, p_\nu \, X^\nu\Bigr).
\end{align*}
So the above lemma expresses exterior powers of clean cycles in terms of convolution products: 
\begin{align*}
 \lambda^{2}(\Lambda) &\;=\; \tfrac{1}{2}\, \bigl( \, \Lambda_{[1,1]} - \Lambda_{[2]} \, \bigr) \\
 \lambda^{3}(\Lambda) &\;=\; \tfrac{1}{6}\, \bigl( \, \Lambda_{[1,1,1]} - 3 \Lambda_{[2,1]} + 2 \Lambda_{[3]} \, \bigr) \\
 \lambda^{4}(\Lambda) &\;=\; \tfrac{1}{24}\, \bigl( \, \Lambda_{[1,1,1,1]} - 6\Lambda_{[2,1,1]} + 3\Lambda_{[2,2]} + 8\Lambda_{[3,1]} - 6\Lambda_{[4]} \, \bigr) \\
  & \;\;\vdots
\end{align*}
\end{ex}

\section{Application to $\scrD$-modules} \label{sec:tannakian-groups}

We now apply the above to characteristic cycles of holonomic $\scrD_A$-modules along the lines of~\cite{KraemerMicrolocal} to get a correspondence between such cycles and Weyl group orbits of weights in representation theory that allows for arbitrary multiplicities.

\subsection{Fiber functors} \label{sec:tannakian}

A module $\scrM \in \Mhol(\scrD_A)$ satisfies $\chi(A, \DR(\scrM))=0$ iff it is negligible in the sense that each of its simple subquotients is stable under translations by a non-zero abelian subvariety~\cite{KrWVanishing, SchnellHolonomic}. Let $\rmS(A) \subset \Mhol(\scrD_A)$ be the Serre subcategory of negligible modules. The quotient category
\[
\rmM(A) \;=\; \Mhol(\scrD_A)/\rmS(A)
\]
is a rigid abelian tensor category for the convolution product 
$*$ induced by the group law $A\times A \to A$~\cite{KraemerMicrolocal}. The discussion in the previous section is motivated by 

\begin{thm}  \label{thm:microlocalization}
Let $\eta \in V$ denote the generic point. Then there exists a $\bbC$-linear exact tensor functor
\[
 \mu: \quad \rmM(A) \;\longrightarrow\; \VB(A, \eta)\smallskip 
\]
such that $\CC(\scrM) \equiv \cleansupp(\mu(\scrM))$ modulo negligible cycles for all $\scrM \in \rmM(A)$. 
\end{thm}

{\em Proof.} Microlocalization even allows to define a tensor functor to an analogous category of local systems~\cite{KraemerMicrolocal}. We will not use the monodromy operation in what follows, so we replace local systems by their sheaves of holomorphic sections. \qed

\medskip 

One reason for working with germs of vector bundles rather than local systems is that there are other interesting tensor functors to which the following discussion applies. So for the rest of this section we fix an {\em arbitrary} $\bbC$-linear exact tensor functor
$F: \rmM(A) \rightarrow \VB(A, \eta)$ where $\eta \in V$ is the generic point; for the moment the specific choice will not matter, although in our applications starting from section~\ref{sec:schottky} we will always use the functor $F=\mu$ from theorem~\ref{thm:microlocalization}.
Fix $\scrM \in \Mhol(\scrD_A)$, and let $\langle \scrM \rangle \subset \rmM(A)$ be the smallest rigid abelian tensor subcategory containing~$\scrM$. It can be described explicitly as the full subcategory of all subquotients of convolution powers of $\scrM\oplus \scrM^\vee \in \rmM(A)$, and as such
has at most countably many isomorphism classes of objects. So for a very general closed point $u\in V(\bbC)$ the restriction of $F$ to this subcategory factors over a $\bbC$-linear exact tensor functor
\[
F_u: \quad \langle \scrM \rangle \;\longrightarrow\; \VB(A, u).
\]
Taking fibers over $u$ we then obtain a $\bbC$-linear exact tensor functor $\omega_u$ such that the diagram
\[
\xymatrix@M=0.5em@R=2em@C=2em{
 \langle \scrM \rangle \ar[rr]^-{\omega_u} \ar[dr]_-{F_u} && \Vect(\bbC) \\
 & \VB(A, u) \ar[ur]_-{\textnormal{fiber}\atop\textnormal{over $u$}} &  
}
\]
commutes. Thus we recover the following Tannakian description~\cite{KrWVanishing,KraemerMicrolocal}:

\begin{cor} \label{cor:tannakian}
There is a linear algebraic group $G=G(\scrM, u)$ such that $\omega_u$ underlies an equivalence 
\[
 \omega_u: \quad \langle \scrM \rangle \;\stackrel{\sim}{\longrightarrow}\; \Rep(G)
\]
with the tensor category of finite-dimensional complex linear representations of $G$.
\end{cor}

{\em Proof.} The functor $\omega_u$ is faithful since any exact tensor functor between rigid abelian categories with $\End(\one)=\bbC$ is so~\cite[prop.~1.19]{DM}. By th.~2.11 in loc.~cit.~the group $G(\scrM, u)=\Aut^\otimes(\omega_u)$ of its tensor automorphisms is then represented by a linear algebraic group which does the job. \qed

\medskip

Since we work over an algebraically closed field, any two fiber functors on our category are non-canonically isomorphic, and so are the corresponding Tannakian groups. We simply write $G(\scrM)=G(\scrM, u)$ and $\omega = \omega_u$ when the specific choice of the fiber functor is either clear from the context or does not matter.

\subsection{Characteristic cycles and weights}

We want to understand the relation between the above Tannakian groups and the clean cycles that arise as supports of objects in the image of our chosen functor $F: \rmM(A) \longrightarrow \VB(A, \eta)$:

\begin{defn} \label{def:cc}
By the {\em clean characteristic cycle} of $\scrM\in \rmM(A)$ we mean the cycle 
\[ \mathrm{cc}(\scrM) \;=\;  \cleansupp(F(\scrM)) \;\in\; \scrZ(A). \]
\end{defn}

\noindent 
If $F=\mu$ is chosen as in theorem~\ref{thm:microlocalization}, then this definition agrees with the usual characteristic cycle modulo negligible cycles, but unlike the latter it is well-defined not only on $\Mhol(\scrD_A)$ but also on the quotient category $\rmM(A)$; we use the lowercase notation to distinguish our clean characteristic cycle from the usual one. For other choices of~$F$ there is no obvious relation between the two but $\cc$ still behaves like a characteristic cycle, for instance Kashiwara's index formula still holds:

\begin{cor} \label{cor:RG-to-ZA}
For $G=G(\scrM, u)$ as in corollary~\ref{cor:tannakian}, we obtain a $\lambda$-ring homomorphism
\[
\cc: \quad R(G) \;=\; K^0(\Rep(G))\;\longrightarrow\; \langle \cc(\scrM) \rangle \;\subset\; \scrZ(A),
\quad \omega_u(\scrN) \;\mapsto\; \cc(\scrN)
\]
with
\[ 
\dim (\omega_u(\scrN)) \;=\; \deg(\cc(\scrN)) \;=\; \chi(A, \DR(\scrN))
\quad 
\textnormal{\em for all $\scrN \in \langle \scrM \rangle$}. \]
\end{cor}

{\em Proof.} By assumption we have a tensor functor
$F_u: \langle \scrM \rangle \longrightarrow \VB(A, u)$. Any such functor induces a $\lambda$-ring homomorphism between the Grothendieck rings, so the first claim holds by proposition~\ref{prop:lambda}. The formula $\dim(\omega_u(\scrN))=\deg(\cc(\scrN))$ follows from our construction of the fiber functor. On the other hand, since any two fiber functors are isomorphic, we have $\dim(\omega_u(\scrN))=\chi(A, \DR(\scrN))$ by~\cite{KrWVanishing}. \qed

\medskip

This gives a dictionary between characteristic cycles and representation-theoretic weights as follows. With notations as in proposition~\ref{prop:lambda}, take the finitely generated subgroup
\[
 \Gamma \;=\; \Gamma(\Lambda, u) \;\subset\; A(\bbC) \quad \textnormal{for} \quad \Lambda \;=\; \cc(\scrM).
\]
Its torsion part contains the subgroup
\[
 K \;=\; K(\scrM, u) \;=\; \bigl\{ \, a\in \Gamma_\mathit{tors}
 \; \mid \; \Supp(\scrN)=\{a\} \; \textnormal{for some} \;  \scrN \in \langle \scrM\rangle \, \bigr\}.
\]
By~\cite[th.~1.3]{KraemerMicrolocal} this subgroup is Cartier dual to the group of components $G/G^\circ$ of the Tannakian group $G=G(\scrM, u)$, and the direct image of $\scrD$-modules under the isogenies in the diagram
\[
\xymatrix{
 A \ar[dr]_-f \ar[rr]^{h\;=\;g\circ f} &&
 A/\Gamma_\mathit{tors} \\
 & A/K \ar[ur]_-g &
}
\]
induces isomorphisms
$G(h_*\scrM, u) \stackrel{\sim}{\longrightarrow} G(f_*\scrM, u) \stackrel{\sim}{\longrightarrow} G^\circ \subseteq G$ onto the connected component of our Tannakian group.
This is particularly useful in the reductive case because representations of connected reductive groups are determined completely by their weights with respect to a maximal torus (example~\ref{ex:character-ring}). In any case we obtain:

\begin{thm} \label{thm:weyl-orbits}
In the above setting, for any maximal torus $T \hookrightarrow G = G(\scrM, u)$ one has an epimorphism 
$ 
 p: X = \Hom(T, \bbG_m)  \twoheadrightarrow \Gamma_\mathit{free} 
$
such that the following diagram commutes:
\[
\xymatrix@M=0.5em@R=2em@C=3em{
R(G) \ar[r]^-{(-)|_{G^\circ}} \ar[d]_-{\cc}
& R(G^\circ) \ar[r]^-\sim \ar[d]_-{\cc}
& R(G^\circ) \ar[r] \ar[d]_-{\cc}  
&\bbZ[X] \ar@{..>}[d]^-p \\
\langle \cc \scrM \rangle \ar[d]_-{(-)_u} \ar@{->>}[r]^-{f_*} 
& \langle f_*\cc\scrM \rangle \ar@{->>}[r]^-{g_*}
& \langle h_*\cc\scrM \rangle \ar[r]^-{(-)_u} 
& \bbZ[\Gamma_\mathit{free}] \ar@{=}[d] \\
\bbZ[\Gamma] \ar@{->>}[rrr]^-{h_*} 
&&& \bbZ[\Gamma_\mathit{free}]
}
\]
\end{thm}

{\em Proof.} 
By the above discussion $f_*\circ \cc$ and $h_*\circ \cc$ factor over $R(G^\circ)$ as in the two upper left squares of the diagram. Furthermore, by construction we have a fiber functor
$
 \omega_u: \Rep(G) \rightarrow \Vect_\Gamma(\bbC)
$
with values in $\Gamma$-graded vector spaces.~The latter are naturally representations of the Cartier dual $\Hom(\Gamma, \bbG_m)$, which is a subgroup of multiplicative type. Its connected component is a subtorus
$
 \Hom(\Gamma_\mathit{free}, \bbG_m) \hookrightarrow G
$
and up to conjugacy we may assume that it sits in a given maximal torus. We then get the desired epimorphism $X\twoheadrightarrow \Gamma_\mathit{free}$ on character groups. \qed

\medskip 

If $\scrM \in \Mhol(\scrD_A)$ is semisimple, then the group $G=G(\scrM, u)$ is reductive. In this case, let us denote by
$
 W(G) = N_{G^\circ}(T)/Z_{G^\circ}(T)
$
the Weyl group of its connected component. The above dictionary then takes the following form:

\begin{cor}
In the above setting, any orbit $O\subseteq X$ of the Weyl group $W(G)$ is realized geometrically by a unique clean effective conic Lagrangian cycle in the sense that
$p(O) = \Lambda_u$ for a unique $\Lambda \in \scrL(A)$.
\end{cor}

{\em Proof.} By the theory of connected reductive groups, the image of $R(G^\circ)\hookrightarrow \bbZ[X]$ are the Weyl group invariants (example~\ref{ex:character-ring}). \qed

\subsection{An inverse Galois problem} \label{sec:inversegalois}

Applying tensor constructions like exterior powers to a faithful representation will usually destroy the faithfulness. In good cases the outcome will still be a faithful representation of an isogenous group, e.g.~$\Alt^n(\bbC^{2n})$ is a faithful representation of the quotient
$\Sl_{2n}(\bbC)/\mu_n$
by the \mbox{$n$-th} roots of unity. In this case we can recover the original group as the universal cover. It is a natural to ask if this abstract Lie-theoretic procedure has a geometric incarnation:
Let $\scrM \in \Mhol(\scrD_A)$. If $G=G(\scrM, u)$ is connected and semisimple, is its universal cover realized as
\[
 \widetilde{G} \;=\; G(\Mtilde, u)\; \twoheadrightarrow\; G
 \quad 
 \textnormal{for some $\Mtilde \in \Mhol(\scrD_A)$ with $\scrM \in \langle \Mtilde \rangle$?}
\]
In this generality this is not true, one can find examples on elliptic curves where one first needs to replace~$\scrM$ by its direct image under an isogeny $A\twoheadrightarrow B$. On the level of characteristic cycles this is the only obstruction: 

\begin{thm} \label{thm:inverse-galois}
Let $\scrM \in \Mhol(\scrD_A)$, and let $p: \widehat{G} \twoheadrightarrow G=G(\scrM, u)$ be an isogeny from an abstract connected reductive group onto its Tannakian group.
Let $S$ be a tensor construction in $r$ variables. If $\scrN \in \langle \scrM \rangle$ and $U_1, \dots, U_r\in \Rep(\widehat{G})$ are such that
\[
 p^*(\omega_u(\scrN)) \;\simeq\;  S(U_1, \dots, U_r), \smallskip
\]
then this isomorphism is realized geometrically by cycles $\Lambda_i\in \langle \cc(\scrM)\rangle \subset \scrZ(A)$ with
\[
  [e]_* \cc(\scrN) \;=\; S(\Lambda_1, \dots, \Lambda_r)
  \quad \textnormal{\em for} \quad
  e \;=\; \deg(p).
\]
\end{thm}

{\em Proof.} For any isogeny $\widehat{G} \twoheadrightarrow G$ the image $T\subset G$ of a maximal torus $\widehat{T} \subset \widehat{G}$ is again a maximal torus, and the degree of the induced isogeny between these tori is again $e$. Then 
$
 X = \CartierDual{T} \hookrightarrow \widehat{X} = \CartierDual{\widehat{T}}
$
is an index $e$ subgroup and hence 
\[
 \bbZ[e\widehat{X}] \;=\;
 \mathrm{im} 
 \Bigl( \!\!
 \xymatrix@M=0.5em{
  \bbZ[\widehat{X}] \ar[r]^-{[e]_*} & \bbZ[\widehat{X}] 
 }
 \!\! \Bigr) \;\subseteq\; \bbZ[X]. \smallskip
\]
Recalling from the proof of lemma~\ref{lem:schur} that on character rings we have $\Psi^e = [e]_*$, we get a factorization
\[
\Psi^e: \quad R(\widehat{G}) \;\longrightarrow\; R(G) \;\subseteq\; R(\widehat{G})
\]
from the diagram
\[
\xymatrix@M=0.6em@C=2em@R=2.5em{
 & R(\widehat{G}) \ar[rrr]^-\sim &&& \bbZ[\widehat{X}]^W \\
 R(\widehat{G}) \ar[rrr]^-\sim \ar[ur]^-{\Psi^e} \ar@{..>}[dr]_-\exists & \ar[u] && \bbZ[\widehat{X}]^W \ar[ur]^-{[e]_*} \ar@{..>}[dr]_-\exists & \\
 & R(G) \ar[rrr]^-\sim \ar@{^{(}-}[u] &&& \bbZ[X]^W \ar@{^{(}->}[uu] 
}
\]
%
%
where $W=W(\widehat{G})=W(G)$ denotes the Weyl group of our two groups. 
Hence we put
\[
 V_i \;:=\; \Psi^e(U_i) \;\in\; R(G).
\]
In general these are only virtual representations, formal $\bbZ$-linear combinations of representations, but this is enough for our purpose. Since the Adams operations are homomorphisms of $\lambda$-rings, the naturality of tensor constructions with respect to such homomorphisms gives \smallskip
\[
 \Psi^e(\omega_u(\scrN)|_{\widehat{G}}) \;=\;
 \Psi^e (S(U_1, \dots, U_r)) \;=\; S(\Psi^e(U_1), \dots, \Psi^e(U_r)) \;=\; S(V_1, \dots, V_r)
\smallskip 
\]
in $R(G)$. Applying the $\lambda$-ring homomorphism $\cc: R(G) \to \langle \scrM \rangle$ from corollary~\ref{cor:RG-to-ZA} we obtain \smallskip 
\[
 [e]_* \cc(\scrN) \;=\; S(\Lambda_1, \dots, \Lambda_r)
 \quad \textnormal{for the cycles} \quad \Lambda_i \;:=\; \cc(V_i) \;\in\; \langle \cc(\scrM) \rangle \smallskip
\]
as desired, where for the left hand side we have used the description of the Adams operations on clean cycles in proposition~\ref{prop:lambda}. Note that we have {\em not} touched the Adams operations on the right hand side: The homomorphism $\cc$ is only defined on the subring $R(G)\subseteq R(\widehat{G})$, so we can apply it to $V_i = \Psi^e(U_i)$ but not to $U_i$.  \qed

\medskip 

\begin{rem} \label{rem:effective}
In the above result the group $G=G(\scrM)$ must be connected. To reduce to the connected case we pass to $A/K$ for the subgroup $K=K(\scrM)\subseteq \Gamma_\mathit{tors}$:
\[
 \xymatrix@M=0.5em{
 A \ar[r]^-f \ar@/_2pc/[rr]|-{\;h\;} & A/K \ar[r]^-g & A/\Gamma_\mathit{tors} 
 } \medskip
\]
The cycles $\Lambda_i$ in the above result need not be effective, but the~$g_* \Lambda_i$ are effective since they come from nonnegative linear combinations of weights in the following diagram:
\[
\xymatrix@M=0.5em{
 R(G) \ar[d]_-{h_*\cc} \ar[r]^-{\mathit{ch}} & \bbZ[X] \ar[d]^-p & \bbN_0[X] \ar@{_{(}->}[l] \ar[d]\\
 \langle h_* \cc(\scrM) \rangle \ar[r]^-{(-)_u} & \bbZ[\Gamma_\mathit{free}] & \bbN_0[\Gamma_\mathit{free}] \ar@{_{(}->}[l]
 }
\]
We refer to this situation by saying that the $\Lambda_i\in \scrL(A)$ are {\em effective up to isogeny}.
\end{rem} 

\medskip

We will apply the above results for the tensor construction $S(V)=\Alt^n(V)$ in theorem~\ref{thm:fake-jacobian} and for $S(V_1, V_2)=V_1\boxtimes V_2$ in proposition~\ref{prop:product}.

\section{Chern-Mather classes} \label{sec:chernmather}
 
For the rest of this paper we will always use the fiber functors from theorem~\ref{thm:microlocalization}, so we work in the subring $\scrL(A)\subset \scrZ(A)$ of clean {\em conic Lagrangian} cycles. This allows to use tools from conormal geometry. We will see that Chern-Mather classes define a ring homomorphism from certain subrings of $\scrL(A)$ to truncations of the Chow ring $\CH_\bullet(A)$ endowed with the Pontryagin product.

\subsection{Generic transversality} 

Recall that any conic Lagrangian subvariety of the cotangent bundle arises as the conormal variety $\Lambda = \Lambda_Z \subset A\times V$ to a closed subvariety $Z\subset A$~\cite[lemma~3]{Kennedy}. We denote its image in the projective cotangent bundle by
\[
 \bbP \Lambda \;\subset\;  A\times \bbP V
\]
and extend the notation additively to conic Lagrangian cycles. Let $p: A\times \bbP V \twoheadrightarrow A$ be the projection.

\begin{defn} 
The {\em Chern-Mather classes} of a conic Lagrangian cycle $\Lambda$ on $T^*A$ are defined by
\[
 c_{M,d}(\Lambda) \;=\; {p_*} \bigl([\bbP \Lambda] \cdot [A\times H_d] \bigr) \;\in\; \CH_d(A) \smallskip
\]
where $H_d \subseteq \bbP V$ is a general linear subspace of dimension $d \in \{0,1,\dots, g-1\}$. 
\end{defn} 

For conormal varieties $\Lambda = \Lambda_Z$ with $Z\neq A$, the triviality of the cotangent bundle implies that these classes coincide with the dual Chern-Mather classes of $Z$ in the sense of~\cite[lemme 1.2.1]{SabbahConormaux} and~\cite[lemma~1]{Kennedy}, except that we consider them as classes on the ambient abelian variety. By Kleiman's generic transversality theorem they are all represented by effective cycles, and they are non-trivial if $\Lambda_Z\subset T^*A$ is clean in the sense of section~\ref{sec:indexformula}:

\begin{lem} \label{lem:mather-transversal}
Let $Z\subset A$ be an irreducible subvariety of dimension $d<g$. \smallskip
\begin{enumerate} 
\item 
For any open subset $U\subset Z^{reg}$ the intersection $\bbP \Lambda_{Z} \cap (U\times H_i)$ is transversal or empty, so 
\[
 c_{M,i}(\Lambda_Z) \;=\; p_* \bigl[ \, \overline{\bbP \Lambda_{Z} \cap (U\times H_i)} \, \bigr]. 
\]
\item 
In particular $d = \max\bigl\{i \mid c_{M, i}(\Lambda_Z)\neq 0 \bigr\}$ and $c_{M, d}(\Lambda_Z) = [Z]$. \medskip 
\item 
If the map $\gamma_Z: \bbP \Lambda_Z \twoheadrightarrow \bbP V$ is dominant, then $c_{M, i}(\Lambda_Z) \neq 0$ for all $i\leq d$. \medskip 
\end{enumerate}
\end{lem}

{\em Proof.} For the first two parts see~\cite[prop.~2.8]{SchuermannTibar}~\cite{KleimanTransversality}. If the map $\gamma_Z$ is dominant, then 
\[ 
 W_i \;=\; \bbP \Lambda_Z \cap (A\times H_i) 
 \quad \textnormal{is nonempty for all $i < g$}.
\]
But we know from the second part that the projection
$
 p: W_i \twoheadrightarrow p(W_i) \subset A
$
is generically finite over its image for $i=d$, and then the same holds for all $i\leq d$ because the $H_i\subset \bbP V$ have been chosen generically. \qed

\medskip

\subsection{Vogel's intersection algorithm}
In order to compute Chern-Mather classes for divisors one can apply Vogel's intersection algorithm~\cite{Vogel,FlennerJoins,FlennerCarrollVogel}. Let us briefly recall how this works in the special case required here: Consider an irreducible projective variety $Z$ of dimension $n$. Fix $\scrL\in \Pic(Z)$, and let $W\subseteq H^0(Z, \scrL)$ be a subspace such that the map
\[
 |W|: \quad Z \;\dashrightarrow\; \bbP W^*
\]
is generically finite and dominant. Then inductively for $n\geq i\geq 0$ we get effective cycles $V_i$ and $R_i$ of pure dimension $i$ such that \smallskip
\begin{itemize}
\item $V_n = \varnothing$ and $R_n = [Z]$, \smallskip
\item $V_i + R_i = R_{i+1} \cap \mathrm{div}(s_i)$ for the divisor of a generic $s_i \in W \subset H^0(Z, \scrL)$, \smallskip
\item $\Supp(V_i)$ is entirely contained in the base locus of $|W|$, \smallskip
\item $\Supp(R_i)$ does not have any component contained in this base locus. \medskip
\end{itemize}
The $V_i$ are called {\em Vogel cycles} and the $R_i$ are the {\em residual Vogel cycles}. Note that their classes
\[
 v_i(\scrL, W) \;=\; [V_i], \quad
 r_i(\scrL, W) \;=\; [R_i] \;\in\; \CH_i(Z)
\]
modulo rational equivalence do not depend on the choice of the generic sections~$s_i$ in the construction. If $Z\subset A$ is a subvariety of our abelian variety, we use the same notation for the images of these cycles in $\CH_i(A)$. They can be expressed via Segre classes but have the advantage of being represented by {\em effective} cycles. As in~\cite{CGS} we have

\begin{cor} \label{cor:chern-gauss}
Let $\scrL = \scrO_A(Z)|_Z$ where $Z\subset A$ is an irreducible reduced ample divisor, and denote by
\[
 W \;=\; \bigl\langle \partial_\nu \vartheta \mid \nu = 1,2,\dots, g \bigr\rangle \;\subseteq\; H^0(Z, \scrL) \medskip
\]
the subspace spanned by the derivatives of a section $\vartheta \in H^0(A, \scrO_A(Z))$ with zero locus $\mathrm{div}(\vartheta)=Z$. Then
\[ c_{M,i}(\Lambda_Z) \;=\; 
\begin{cases} 
\; r_i(\scrL, W) & \textnormal{\em for all $i\in\{0,1,\dots, g-1\}$}, \\
\; [Z]^{g-1} & \textnormal{\em for all $i>\dim \Sing(Z)$}.
\end{cases}
\]
\end{cor} 

{\em Proof.} For a divisor the projection $p: \bbP\Lambda_Z \to Z$ is a birational map, and the composite
\[
 \gamma_Z: \quad Z \;\dashrightarrow \; \bbP\Lambda_Z \; \hookrightarrow \; A\times \bbP V \;\twoheadrightarrow \; \bbP V
\]
is the Gauss map sending a smooth point of the divisor to the conormal direction at that point. Up to a scalar there is a unique $\vartheta \in H^0(A, \scrO_A(Z))$ with $div(\vartheta)=Z$, and 
$\gamma_Z(z) = \bigl[ \partial_1 \vartheta(z) : \cdots : \partial_g \vartheta(z)\bigr]$
for any smooth point of the divisor. So the Gauss map is the rational map defined by the linear series $|W|$, its base locus is the singular locus of the divisor, and for an ample divisor it is dominant and generically finite~\cite[sect.~4.4]{BL_ComplexAbelian}. Now apply part (1) of the previous lemma. \qed

\medskip 

\subsection{The Pontryagin product}

Let us now consider Chern-Mather classes of convolutions.
For a conic Lagrangian subvariety $\Lambda \subset T^* A$ its total Chern-Mather class is
\[
 c_M(\Lambda) \;=\; c_{M,0}(\Lambda) + c_{M,1}(\Lambda) + \cdots + c_{M, {g-1}}(\Lambda) \;\in\; \CH_\bullet(A).
\]
By the above this class is given by an effective cycle, which vanishes iff $\Lambda \subset T^* A$ is the zero section. So if $A$ has nontrivial proper abelian subvarieties, then there are negligible conic Lagrangian subvarieties whose total Chern-Mather class does not vanish. In what follows we restrict ourselves to clean cycles and consider the group homomorphism
\[
 c_M: \quad \scrL(A) \;\longrightarrow\; \CH_\bullet(A)
\]
which sends a clean cycle to its total Chern-Mather class. Ideally we would like this to be a ring homomorphism when the target is endowed with the Pontryagin product
\[
 *: \quad 
 \CH_{a}(A) \times \CH_{b}(A) \; \longrightarrow \; \CH_{a+b}(A), \quad \alpha * \beta \;=\; a_*(\alpha\boxtimes \beta), 
\]
but this cannot be true: In the definition of the convolution product $\circ$ on~$\scrL(A)$ we have only worked with germs of subvarieties, whereas in general the Pontryagin product of clean cycles may involve negligible contributions. To control the latter we need estimates on positive-dimensional fibers of Gauss maps. For $1\leq d \leq g-1$, let us denote by
\[
 \scrL^{>d}(A) \;\subseteq\; \scrL(A)
\]
the subgroup generated by all clean conic Lagrangian subvarieties $\Lambda \subset T^* A$ whose Gauss map
$
 \gamma_{\bbP \Lambda}:  \bbP \Lambda \rightarrow \bbP V 
$
restricts to a finite morphism over the complement of a proper closed subset $S\subset \bbP V$ of codimension $\codim(S, \bbP V) > d$. Since the class of finite morphisms is stable under base change, one easily checks that this defines a filtration by subrings
\[
 \scrL(A) \;=\; \scrL^{>1}(A) \;\supseteq\; \scrL^{>2}(A) \;\supseteq\; \cdots \;\supseteq\; \scrL^{>g-1}(A) 
\]
where the final step on the right is generated as a group by the conic Lagrangian subvarieties whose Gauss map is finite. Similarly, for any $d$ consider the quotient ring
\[
 \CH_{\leq d}(A) \;=\; \CH_\bullet(A)/\CH_{>d}(A)
 \quad \textnormal{by the ideal} \quad
 \CH_{>d}(A) \;=\; \bigoplus_{i>d} \CH_{i}(A), 
\]
then we obtain:

\begin{lem} \label{lem:mather-ringhomo}
For any $d\in \{1,2,\dots, g-1\}$ the Chern-Mather class gives rise to a ring homomorphism
\[ c_M: \quad \scrL^{>d}(A) \; \longrightarrow\; \CH_{\leq d}(A). \]
More precisely, we have
\[
 c_M(\Lambda_1 \circ \Lambda_2) - c_M(\Lambda_1)*c_M(\Lambda_2)
 \in \CH_{>d}(A)
 \;\;\;\textnormal{\em for all} \;\;\; \Lambda_1 \in \scrL^{>d}(A), \; \Lambda_2\in \scrL(A).
\]
\end{lem}

{\em Proof.} The first claim follows from the second. For the latter we may assume by additivity that $\Lambda_1, \Lambda_2 \subset T^* V$ are irreducible subvarieties. The main point will then be to give an intersection-theoretic expression for the class of $\Lambda_1 \circ \Lambda_2$ in~$\CH_{g-1}(A\times \bbP V)$. The diagonal  
\[
 \delta: \quad \Delta \;=\; A^2 \times \bbP V \;\hookrightarrow\;  A^2\times \bbP V \times \bbP V
\]
is a regular embedding of codimension $g-1$, so by~\cite[lemma~7.1]{Fulton} any irreducible component 
\[
 W \; \subset \; \delta^{-1}\bigl( \bbP \Lambda_1 \times \bbP \Lambda_2 \bigr) \;=\; \bbP \Lambda_1 \times_{\bbP V} \bbP \Lambda_2  
\]
has $\dim(W) \geq g-1$. Strict inequality can only occur if the image of $W$ in $\bbP V$ is contained in the locus where at least one of the Gauss maps $\bbP \Lambda_i \to \bbP V$ has positive-dimensional fibers. Extending the notions from section~\ref{sec:indexformula} to this case, we say that \smallskip 
\begin{itemize}
\item $W$ is {\em clean} if the projection $W\to \bbP V$ is dominant and hence restricts to a finite \'etale cover over an open dense subset of the target. \smallskip 
\item $W$ is {\em negligible} otherwise, in particular whenever $\dim W > g-1$. \smallskip
\end{itemize}
Now $\Lambda_1 \circ \Lambda_2$ is a linear combination of clean conic Lagrangian subvarieties. After projectivization we obtain by definition of the product $\circ$ that $\bbP(\Lambda_1\circ \Lambda_2)=\varpi_*(\alpha)$, where
\[
 \varpi: \quad A^2 \times \bbP V \;\longrightarrow\; A\times \bbP V
\]
is the addition map and
\[
 \alpha \;=\; \sum_{W\mathit{clean}} [W] \;\in\; \CH_{g-1}(A^2\times \bbP V)
\]
is the sum of all clean irreducible components $W\subset \delta^{-1}(\bbP \Lambda_1 \times \bbP \Lambda_2)$. Note that each of these components is reduced and enters with multiplicity one because the projection $W\to \bbP V$ is finite \'etale over an open dense subset of the target. On the other hand the class
\[
 \beta \;=\; \delta^!([\bbP \Lambda_1 \times \bbP \Lambda_2]) \;\in\; \CH_{g-1}(A^2\times \bbP V)
\]
is represented by a cycle which might contain also negligible summands, but it follows from basic intersection theory that $\alpha$ and $\beta$ agree on all clean components. So we may write
\[
 \alpha - \beta \;\;=\!\!\! \sum_{W \mathit{negligible}} m_W \cdot \beta_W
 \quad \textnormal{with} \quad 
 \beta_W \;\in\; \mathrm{Im}\Bigl(\CH_{g-1}(W)\to \CH_{g-1}(A^2 \times \bbP V)\Bigr)
\]
for certain multiplicities $m_W\in \bbZ$.

\medskip 

Now let $S\subset \bbP V$ be a closed subset such that the Gauss map $\gamma: \bbP \Lambda_1 \to \bbP V$ restricts to a finite morphism over $U=\bbP V \setminus S$. Since finite morphisms are stable under base change,
\[
 \gamma \times \id: \quad \delta^{-1}(\bbP \Lambda_1 \times \bbP \Lambda_2)|_U \; \longrightarrow\; \bbP \Lambda_2|_U
\]
is then also a finite morphism. Since the target is irreducible of dimension $g-1$, it follows that 
\[
 (\gamma \times \id)(W) \;\subseteq\; S
 \quad \textnormal{for any negligible component} \quad 
 W \;\subset\; \delta^{-1}(\bbP \Lambda_1 \times \bbP \Lambda_2). \medskip
\]
If $\codim(S, \bbP V)>d$, it follows that via the K\"unneth decomposition for the Chow ring of a projective bundle \medskip
\[
 \alpha - \beta \;\in\;
 \bigoplus_{i>d} \; \CH_{i}(A^2) \otimes \CH_{g-1-i}(\bbP V)
 \;\subset\; \CH_{g-1}(A^2\times \bbP V). 
\]
To conclude the proof, one easily checks from the definition of the classes $\alpha$ and $\beta$ that
\begin{align*}
 \varpi_*(\alpha) &= \sum_i c_{M,i}(\Lambda_1\circ \Lambda_2) \otimes H^i \\
 \varpi_*(\beta) &= \sum_{i,j} (c_{M,i}(\Lambda_1)*c_{M,j}(\Lambda_2))\otimes H^{i+j}
\end{align*}
where $H\in \CH_{g-2}(\bbP V)$ denotes the class of a hyperplane. \qed

\medskip

\begin{cor} \label{cor:mather-ringhomo}
If $A$ is a simple abelian variety, the total Chern-Mather class gives a ring homomorphism 
\[ c_M: \quad \scrL(A) \;\longrightarrow\; \CH_{\bullet}(A)/\CH_g(A). \]
\end{cor}

{\em Proof.} If there exists a subvariety $Z\subset A$ whose Gauss map $\gamma_Z: \bbP \Lambda_Z \twoheadrightarrow \bbP V$ is dominant but admits a positive-dimensional fiber $\gamma_Z^{-1}(\xi)$, then this fiber generates a nontrivial abelian subvariety. On a simple abelian variety this is impossible, so we get $\scrL(A) = \scrL^{>g-1}(A)$ and the previous lemma applies. \qed

\section{The Schottky problem} \label{sec:schottky}

One motivation for the inverse Galois problem from section~\ref{sec:inversegalois} is the Schottky problem whether the Tannakian formalism detect Jacobians among all ppav's. In this section we discuss a possible approach to this question combining theorem~\ref{thm:inverse-galois} with the results of the previous section about Chern-Mather classes.

\subsection{Tannakian Schottky} 

Let $A$ be a ppav with theta divisor $\Theta \subset A$. The theta divisor is determined by the polarization only up to a translation, and in what follows we fix one of the $2^{2g}$ symmetric translates; in section~\ref{sec:translates} we will see that the specific choice does not matter too much. If $\Theta$ is smooth, we know by~\cite{KraemerMicrolocal} that
\[
 G(\delta_\Theta, u) \;\simeq\;
 \begin{cases}
 \Sp_{g!}(\bbC) & \textnormal{if $g$ is even}, \\
 \SO_{g!}(\bbC) & \textnormal{if $g$ is odd},
 \end{cases}
\]
and $\omega_u(\delta_\Theta)$ is the standard representation. For theta divisors on special ppav's the group can be much smaller:

\begin{ex} If $(A, \Theta)=Jac(C)$ is the Jacobian of a smooth projective curve~$C$ of genus $g=n+1$, we have the resolution 
$C^{(n)} \twoheadrightarrow \Theta$. By the decomposition theorem we get an inclusion
\[ 
 \delta_\Theta \;\hookrightarrow\;  
 \Alt^{* n}(\delta_C) 
\]
and by~\cite[sect.~6]{KrWSmall}~\cite{WeBN} it follows that we have an isogeny $\tilde{G} = G(\delta_C, u)\twoheadrightarrow G(\delta_\Theta, u)$ with 
\begin{equation} \tag{$\star$} \label{eq:jacobian}
 \tilde{G} \;\simeq\; 
 \begin{cases} 
 \Sl_{2n}(\bbC) \\
 \Sp_{2n}(\bbC)
 \end{cases}
 \textnormal{acting via} \quad
 \omega_u(\delta_\Theta)|_{\tilde{G}} \;\simeq\;
 \begin{cases}
 \Alt^n(\bbC^{2n}) \\
 \Alt^n(\bbC^{2n})/\Alt^{n-2}(\bbC^{2n}) 
 \end{cases}
\end{equation}
in the nonhyperelliptic resp.~hyperelliptic case.
\end{ex} 

It is natural to ask whether this characterizes Jacobians. We call a ppav $(A, \Theta)$ of dimension $g=n+1$ a nonhyperelliptic or hyperelliptic {\em fake Jacobian} if for some symmetric translate of the theta divisor the semisimple group $G(\delta_\Theta, u)$ is connected and its universal cover acting on $\omega_u(\delta_\Theta)$ has the form~(\ref{eq:jacobian}). In the moduli space $\scrA_g$ consider the loci 
\[
 \scrJ_g \;\subseteq\; \scrJ_{g,\mathit{fake}} 
\]
of Jacobians of smooth projective curves and of fake Jacobians; the latter is a locally closed algebraic subset by~\cite[prop.~7.4]{KrWSchottky}. We would like to see if the two loci are the same, so for any fake Jacobian we need to find  a candidate for the curve whose Jacobian it should be. A natural guess is provided by

\begin{thm} \label{thm:fake-jacobian}
If $(A, \Theta)$ is a fake Jacobian, then there exists a cycle $\Lambda \in \langle \cc(\delta_\Theta) \rangle$, effective up to isogeny, such that \smallskip
\[ 
  [e]_*\cc(\delta_\Theta) \;=\; 
  \begin{cases} 
  \mathrm{Alt}^{*(g-1)}(\Lambda) & \textnormal{\em with $e\;=\;g-1$} \\ 
  \mathrm{Alt}^{*(g-1)}(\Lambda) - \mathrm{Alt}^{*(g-3)}(\Lambda) & \textnormal{\em with $e\;=\;\gcd(2, g-1) \in \{1,2\}$} 
  \end{cases} 
\]
in the nonhyperelliptic resp.~hyperelliptic case.
\end{thm}

{\em Proof.} Apply theorem~\ref{thm:inverse-galois} to the connected reductive group $G=G(\delta_\Theta, u)$. We have $f=\id_A$ and by direct inspection the degree $e$ of the universal covering map is the given one. So the claim follows by taking $S(U)=\Alt^{g-1}(U)$ respectively $S(U)=\Alt^{g-1}(U)/\Alt^{g-3}(U)$; remark~\ref{rem:effective} says that the resulting cycle is effective up to isogeny. \qed

\medskip

If $(A, \Theta)$ is the Jacobian of a smooth projective curve, then the above cycle $\Lambda$ is indeed the conormal variety to the image of the curve under the Abel-Jacobi map and this gives a constructive proof of Torelli's theorem~\cite{KraemerMicrolocal}. For fake Jacobians the situation is less clear. We would be done if we could show that the image $Z\subset A$ of~$\Lambda\subset A\times V$
is a curve, as then   
$ \Theta = Z+\cdots + Z $ 
is a sum of $g-1$ copies of this curve and hence a Jacobian by~\cite{SchreiederCurveSummands}. Unfortunately it is hard to control convolutions with excess dimension: How to show $\dim(Z)=1$?

\subsection{A computation of Chern-Mather classes} \label{sec:weak-schottky}

One way to control dimensions is to look at Chern-Mather classes. Let us illustrate this technique by a simple example. For any $g$ the locus of Jacobians $\scrJ_g$ is an irreducible component of the Andreotti-Mayer locus \medskip
\[
 \scrN_{g-4} \;=\; \{ (A, \Theta) \in \scrA_g \mid \dim \Sing(\Theta) \geq g-4 \},  \medskip
\]
see~\cite{AM}. So for abelian fivefolds we get 

\begin{cor} \label{cor:schottky}
For $g=5$, any nonhyperelliptic fake Jacobian lies in $\scrN_{g-4}$ and so the closure of the locus of fake Jacobians in $\scrA_g$ contains the closure of the Jacobian locus as an irreducible component.
\end{cor}

{\em Proof.} Let $(A, \Theta)$ be a nonhyperelliptic fake Jacobian of genus $g=5$. In what follows we work up to numerical equivalence and therefore view Chern-Mather classes $c_M(-)\in H_{2\bullet}(A, \bbZ)$ as classes in the even homology ring without further notice. This being said, let us write $\cc(\delta_\Theta) = \Lambda_\Theta + \Lambda'$ where $\Lambda'$ is supported over the singular locus of the theta divisor. If the latter were empty or finite, we would get that
\[
 c_{M,1}(\cc(\delta_\Theta)) \;=\; c_{M,1}(\Lambda_\Theta) \;=\; [\Theta]^4
 \quad \textnormal{in} \quad 
 H_2(A, \bbZ),
\]
where the second equality again uses our assumption on isolated singularities and corollary~\ref{cor:chern-gauss}. For $\Lambda \in \scrL(A)$ as in theorem~\ref{thm:fake-jacobian}, we compute the Chern-Mather classes of
\[
 [g-1]_* \cc(\delta_\Theta) \;=\; \Alt^{*(g-1)}(\Lambda) 
\]
in two different ways, looking at both sides of the above equation. For the right hand side put $c_i = c_{M,i}(\Lambda)$. By construction $c_0 = 2g-2 = 8$. Now consider the cycles
\[
 \Lambda_{[\beta]} \;=\; [\beta_1]_* (\Lambda) \circ \cdots \circ [\beta_\ell]_* (\Lambda)
 \quad \textnormal{for partitions} \quad 
 \beta \;=\; (\beta_1, \dots, \beta_\ell).
\]
Lemma~\ref{lem:mather-ringhomo} allows to easily compute their Chern-Mather classes in degree $d=1$, yielding
\begin{align*} 
 c_{M,1}(\Lambda_{[1,1,1,1]} ) &\;=\; 4\cdot c_0^3 \cdot c_1 \;=\; 2048 \cdot c_1, \\
 c_{M,1}(\Lambda_{[2,1,1]} ) &\;=\; 6 \cdot c_0^2 \cdot c_1\;=\; 384 \cdot c_1, \\
 c_{M,1}(\Lambda_{[2,2]} ) &\;=\; 8 \cdot c_0 \cdot c_1\;=\; 64\cdot c_1, \\
 c_{M,1}(\Lambda_{[3,1]}) &\;=\; 10 \cdot c_0 \cdot c_1 \;=\; 80\cdot c_1, \\
 c_{M,1}(\Lambda_{[4]}) &\;=\; 16 \cdot c_1\;=\; 16\cdot c_1.  
\end{align*}
By example~\ref{ex:elementary-symmetric} then \medskip
\[
 c_{M,1}(\Alt^4(\Lambda)) = 
 \tfrac{1}{24}\, c_{M,1}\bigl( \, \Lambda_{[1,1,1,1]} - 6\Lambda_{[2,1,1]} + 3\Lambda_{[2,2]} + 8\Lambda_{[3,1]} - 6\Lambda_{[4]} \, \bigr)  
 = 20\cdot c_1. \medskip
\]
On the other hand \medskip 
\[
 c_{M,1}([g-1]_* \cc(\delta_\Theta)) \;=\; (g-1)^2 \cdot c_{M,1}(\Lambda_\Theta) \;=\; 16 \cdot [\Theta]^4 \medskip
\]
by corollary~\ref{cor:chern-gauss}, again using our assumption that the theta divisor has at most isolated singularities.
Altogether then $16 \cdot [\Theta]^4  = 20 \cdot c_1$. But this leads to the contradiction
\[
 c_1 \;=\; \frac{4}{5} \cdot [\Theta]^4 \;\notin\; H_2(A, \bbZ)
\]
since $[\Theta]^4$ is not divisible by five, being $24$ times a primitive class in $H_2(A, \bbZ)$. \qed

\subsection{Dependence on the translate} \label{sec:translates}

Let us briefly explain to what extent the above depends on the translate of the theta divisor. Taking $\Theta \subset A$ to be symmetric determines $G(\delta_\Theta, u)$ except for a small issue about connected components: 

\begin{ex}
For $g=1$ the theta divisor is a $2$-torsion point, so the group~$G(\delta_\Theta)$ will be trivial if we choose the point to be the origin, while otherwise it will have order two. A more interesting example is the intermediate Jacobian of a smooth cubic threefold. Here $g=5$ and the theta divisor has a unique singularity $x\in A[2]$, whence
\[
 G(\delta_\Theta) \;\simeq\;
\begin{cases}
 E_6(\bbC) & \textnormal{for $x=0$}, \\
 E_6(\bbC) \times \bbZ/2\bbZ & \textnormal{for $x\neq 0$},
\end{cases}
\]
as one easily sees in terms of the action of $\CartierDual{\langle x\rangle}\subset G(\delta_\Theta)$ in~\cite[ex.~2.2]{KraemerMicrolocal}.
\end{ex}

In the above examples the groups arising from different translates have the same connected component, which is a semisimple group since by the symmetry of the theta divisor its irreducible faithful representation $\omega_u(\delta_\Theta)$ must be self-dual. More generally, let $t_a: A\to A, x\mapsto x+a$ denote the translation by a point $a\in A(\bbC)$ and put
\[
 \scrM_a \;=\; t_a^*(\scrM) \;=\; \scrM * \delta_a \quad \textnormal{for} \quad \scrM \;\in\; \Mhol(\scrD_A).
\]
If $G=G(\scrM)$ is reductive, let
$
 G' = [G^\circ, G^\circ]
$
be the derived group of its connected component. The following shows that even for nonsymmetric theta divisors the situation does not seriously depend on the chosen translate:

\begin{lem} \label{lem:translate}
If $\scrM$ is semisimple, then $G(\scrM)' \simeq G(\scrM_a)'$ for all $a\in A(\bbC)$.
\end{lem}

{\em Proof.}
Put $\scrN = (\scrM \boxtimes \delta_0) \oplus \delta_{(0,a)} \in \Mhol(\scrD_{A\times A})$. The group law $a: A\times A\to A$ gives rise to the following commutative diagram between tensor categories, which by Tannakian duality translates to a diagram of reductive groups:
\[
\vcenter{\vbox{
\xymatrix@R=2em@C=2em@M=0.5em{
 \langle \scrM \oplus \delta_a \rangle & \langle \scrN \rangle \ar@{->>}[l]_-{a_*} \\
 &  \langle \scrM \rangle \ar[ul] \ar[u]_-{(-)\boxtimes \delta_0}
}}}
\qquad 
\vcenter{\vbox{
\xymatrix@R=2em@C=2em@M=0.5em{
G(\scrM\oplus \delta_a) \ar@{->>}[dr] \ar@{^{(}->}[r] & G(\scrM) \times G(\delta_a) \ar@{->>}[d] \\
&  G(\scrM)
}}}
\]
The functor which sends a reductive group to the derived group of its connected component clearly preserves embeddings, epimorphisms and direct products and it sends the multiplicative group $G(\delta_a)$ to the trivial group, hence from the above we get a diagram
\[
\xymatrix@R=2em@C=2em@M=0.5em{
G(\scrM\oplus \delta_a)' \ar@{->>}[dr] \ar@{^{(}->}[r] & G(\scrM)' \ar@{=}[d] \\
&  G(\scrM)'
}
\]
which shows that the diagonal arrow must be an isomorphism. On the other hand, since $\langle \scrM\oplus \delta_a \rangle =\langle \scrM_a\oplus \delta_a\rangle$, we know that $G(\scrM \oplus \delta_a)' \simeq G(\scrM_a\oplus \delta_a)'$ and hence the claim follows by symmetry. \qed

\medskip

\section{Summands of subvarieties}

Our next application is a criterion to detect nontrivial summands in a given subvariety $Z\subset A$. For this we will consider the adjoint representation, but before doing so we collect some general facts about conormal varieties.

\subsection{Geometric nondegeneracy} \label{sec:sum}

An irreducible subvariety $X\subset A$ is said to be~{\em geometrically nondegenerate} if for any epimorphism $A\twoheadrightarrow B$ to an abelian variety the induced morphism $X \to B$ is either surjective or generically finite onto its image~\cite[II.12]{RanSubvarieties}. This is weaker than being {\em nondegenerate} in Ran's sense~\cite[II.1]{RanSubvarieties} but still forces the stabilizer
\[
 \Stab(X) \;=\; \{ \, a\in A \mid X + a = X \,\}
\]
to be finite. The converse does not hold in general: For example, replacing $A$ by a larger ambient abelian variety will destroy the geometric nondegeneracy but does not affect the stabilizer. By~\cite{WeissauerDegenerate}, the stabilizer $\Stab(X)$ is finite if and only if the Gauss map
\[
 \gamma: \quad \Lambda_X \;\subset\; T^*A \;=\; A\times V \;\twoheadrightarrow\; V
\]
is generically finite, i.e.~iff the conormal variety is not negligible. 

\subsection{Two lemmas about conormal varieties}
The following two lemmas allow to control the convolution product of conormal varieties:

\begin{lem} \label{lem:mindim-of-convolution}
If $Z_1, Z_2\subset A$ are irreducible subvarieties with finite stabilizer, every nonnegligible irreducible component $\Lambda_Y \subseteq  \Supp(\Lambda_{Z_1} \circ \Lambda_{Z_2})$ of their convolution satisfies
\[
 d_Y \;\geq\; |d_{Z_1} - d_{Z_2}|,
\]
and equality implies 
\[ Z_1 \;=\; Y - Z_2 \quad \textnormal{or} \quad Z_2 \;=\; Y - Z_1.\]
\end{lem}

{\em Proof.} If $\Lambda_Y$ is an irreducible component of the convolution, then by definition of the convolution there exists a component $\Lambda$ of $\Lambda_{Z_1} \times_V \Lambda_{Z_2}$ that dominates $Y$ via the sum map:
\[
\xymatrix@R=2em@C=2em@M=0.5em{
 \Lambda \ar@{^{(}->}[r] \ar@{->>}[d] & \Lambda_{Z_1}\times_V \Lambda_{Z_2} \ar[d] & ((z_1, \xi), (z_2, \xi)) \ar@{|->}[d] \\
 Y \ar@{^{(}->}[r] & A & z_1 + z_2
}
\]
Composing the inclusion of this component in the fiber product with the projection on the second factor we get a map
$
\Lambda \rightarrow \Lambda_{Z_2}, 
 ((z_1, \xi), (z_2, \xi)) \mapsto (z_2, \xi)
$
which is dominant because both the source and the target are irreducible and over some open dense subset of $V$ both are finite \'etale. In particular, for general $z_2\in Z_2$ we can always find a point $z_1\in Z_1$ with $y=z_1+z_2\in Y$. We therefore have an inclusion~$Z_2 \subseteq Y - Z_1$ and the same holds after interchanging $Z_1$ and $Z_2$. \qed

\begin{lem} \label{lem:multiplicity-one}
If $Z_1, Z_2 \subset A$ are irreducible subvarieties with finite stabilizer and~$\Lambda_Y$ is a component of multiplicity one in the cycle $\Lambda_{Z_1}\circ \Lambda_{Z_2}$, then one has dominant rational maps 
\[
\xymatrix@C=3em@M=0.5em{
 Y \ar@{-->>}[r] & Z_i
}
\quad \textnormal{\em for} \quad i \;=\; 1,2.
\]
\end{lem}

{\em Proof.} As above there exists an irreducible component $\Lambda \subseteq \Lambda_{Z_1}\times_V \Lambda_{Z_2}$ which dominates $\Lambda_Y$. Here
$
 \varpi: \Lambda \twoheadrightarrow \Lambda_Y
$
is generically finite of degree one because of our multiplicity one assumption. Hence $\varpi$ is birational, and composing a rational inverse with the projection on the $i$-th factor we get a rational map $f_i$ as indicated below:
\[
\xymatrix@R=2em@C=1.8em@M=0.5em{
 \Lambda \ar@{^{(}->}[r] \ar[d]_-\varpi & \Lambda_{Z_1} \times_V \Lambda_{Z_2} \ar@{->>}[r] & \Lambda_{Z_i} \ar@{=}[d] \\
 \Lambda_Y \ar@{-->>}[rr]^-{\exists f_i} && \Lambda_{Z_i}
}
\]
As in the previous proof $f_i$ is dominant since it commutes with the projection to~$V$ and both the source and the target are irreducible generically finite covers of the latter. Composing with the map to the abelian variety we get a dominant rational map
$
 \Lambda_Y \dashrightarrow Z_i.
$
Now $\Lambda_Y$ contains the conormal bundle to the smooth locus of~$Y$ as a Zariski open dense subset, and as such it is birational to the product $\tilde{Y} \times \bbC^d$ where $d=g-d_Y$ and $\tilde{Y}\to Y$ is a resolution of singularities. Since any rational map from a smooth variety to an abelian variety extends to a morphism~\cite[th.~3.1]{MilneAV}, we get a morphism $\tilde{Y} \times \bbC^d \rightarrow A$ with image~$Z_i\subset A$. By the universal property of the Albanese the latter factors as
\[
\xymatrix@R=2em@C=2em@M=0.5em{
 \tilde{Y} \times \bbC^d \ar[r] \ar[d]_-{pr} & Alb(\tilde{Y} \times \bbC^d) \ar@{=}[d] \ar[r] & A\\
 \tilde{Y} \ar[r] & Alb(\tilde{Y}) 
}
\]
where $pr$ is the projection. So we get a morphism $\tilde{Y} \rightarrow A$ with image~$Z_i\subset A$.
\qed

\subsection{Summands of subvarieties}

Now let $Z\subset A$ be an irreducible subvariety with trivial stabilizer $\Stab(Z)=\{0\}$ and suppose that there exists a decomposition as a sum
\[ Z \;=\; X+Y \]
of irreducible geometrically nondegenerate subvarieties of positive dimension. Then by~\cite[sect.~8.2]{DebarreComplexTori} the addition morphism $X\times Y \twoheadrightarrow Z$ is generically finite and $Z$ is geometrically nondegenerate. By the decomposition theorem~\cite{BBD} we have an inclusion
\[
 \delta_Z \;\hookrightarrow\; \delta_X * \delta_Y \;\in\; \langle \delta_{XY} \rangle
 \quad \textnormal{where} \quad \delta_{XY} \;:=\; \delta_X \oplus \delta_Y.
\]
So we have a fully faithful embedding as a tensor subcategory $\langle \delta_Z \rangle \hookrightarrow \langle \delta_{XY} \rangle$ stable under subquotients and compatible with our fiber functors. Any such embedding gives an epimorphism
$p: G(\delta_{XY}) \twoheadrightarrow  G(\delta_Z)$ 
\cite[cor.~2.9, prop.~2.21]{DM}, so we may consider
\[ V_S \;=\; \omega(\delta_S) \;\in\; \Rep(G(\delta_{XY})) \quad \textnormal{for} \quad S \;=\; X, Y, Z. \]

\begin{lem}
The above three representations remain irreducible when restricted to the connected component $G(\delta_{XY})^\circ \subseteq G(\delta_{XY})$, and replacing $X, Y, Z$ by translates we may assume 
\[
 G_S \;:=\; G(\delta_S)^\circ 
 \quad \textnormal{\em is semisimple for} \quad 
 S \;\in \; \{ X,Y,Z,XY\}.
\]
\end{lem}

{\em Proof.} The triviality of the stabilizer $\Stab(Z)$ implies that $\Stab(X)$ and $\Stab(Y)$ are trivial. Hence $V_X, V_Y, V_Z$ restrict to irreducible representations of the connected component of the respective Tannaka groups by~\cite[cor.~1.4]{KraemerMicrolocal}. All these groups are quotients of $G(\delta_{XY})$, hence the connected component of the latter also acts irreducibly on the three representations. 

\medskip 

For semisimplicity, recall that any one-dimensional representation of~$G(\delta_{XY})$ is given by a skyscraper sheaf $\delta_a$ on a point $a\in A(\bbC)$~\cite[sect.~3.c]{KraemerMicrolocal}~\cite[lemma~13]{WeissauerDegeneratePerverse}, so
\[
 \det(V_X) \;=\; \omega(\delta_x) \quad \textnormal{and} \quad \det(V_Y) \;=\; \omega(\delta_y)
 \quad \textnormal{for some $x,y\in A(\bbC)$}.
\]
Replacing $X, Y\subset A$ by suitable translates we may assume $x=y=0$ so that the above two determinants become trivial. But by Schur's lemma the center $Z\subset G_{XY}$ acts via scalars on both $V_X$, $V_Y$, and the triviality of the determinant forces these scalars to be roots of unity. Since the action of $G_{XY}$ on $V_X\oplus V_Y$ is faithful, it follows that $Z$ is finite and so the group $G_{XY}$ is semisimple. Since $G_X, G_Y, G_Z$ are quotients of the group $G_{XY}$, they are then semisimple as well. \qed 

\medskip

By semisimplicity, the universal cover $\pi_S: \tilde{G}_S \twoheadrightarrow G_S$ of each of the above groups is again an algebraic group. 
By the correspondence between simply connected semisimple groups and their Lie algebras there exists a unique embedding $i$ making the diagram 
\[
\xymatrix@M=0.5em@C=3em{
\tilde{G}_Z \ar@{^{(}..>}[r]^-{\exists ! i} \ar[dr]_-{\pi_Z} & \tilde{G}_{XY} \ar@{->>}[d]^-{p\circ \pi_{XY}}  \\ 
& G_Z 
}
\]
commute, and we get 

\begin{prop}
The subgroup $\tilde{G}_Z$ acts nontrivially on both $\omega(\delta_X)$, $\omega(\delta_Y)$.
\end{prop}

{\em Proof.} We have a splitting into $G_1=\mathrm{Im}(i)$ and $G_2 = \ker(p\circ \pi_{XY})$ as indicated below:
\[
\xymatrix@M=0.5em{
 G_1\times G_2 \ar@{->}[r]^-\sim \ar@{->}[d]_-{pr_1}
 & \tilde{G}_{XY} \ar@{->>}[r] \ar@{->}[d]^-{\tilde{p}} & G_{XY} \ar@{->}[d]^-p \\
 G_1 \ar[r]^-\sim & \tilde{G}_Z \ar@{->>}[r] & G_Z
}
\]
So any irreducible representation $U\in \Rep(G_{XY})$ can be written as $U \simeq U_1\boxtimes U_2$ where the $U_i\in \Rep(G_i)$ are irreducible representations of the two factors. In our case
\begin{align*} 
V_X &\;\simeq\; V_{X,1}\boxtimes V_{X,2},
& V_Z &\;\subseteq\; V_{X,1}\otimes V_{Y,1},\\
V_Y &\;\simeq\; V_{Y,1}\boxtimes V_{Y,2},
& \one &\;\subseteq\; V_{X,2}\otimes V_{Y,2},
\end{align*}
where $\one \in \Rep(G_2)$ denotes the trivial representation. Note that the last inclusion is equivalent to
\[ V_{X,2} \;\simeq\; \Hom(V_{Y,2}, \one) \] 
because both these representations are irreducible. For sufficiently divisible $n\in \bbN$ and~$S=X,Y$ theorem~\ref{thm:inverse-galois} and remark~\ref{rem:effective} yield effective cycles $\Lambda_{S,i}\in \langle \cc(\delta_{XY})\rangle$ so that
\begin{align*} 
\cc(\delta_{nX}) &\;=\; \Lambda_{X,1}\circ \Lambda_{X,2},
&\cc(\delta_{nZ}) &\;\subseteq\; \Lambda_{X,1}\circ \Lambda_{Y,1}, \\
\cc(\delta_{nY}) &\;=\; \Lambda_{Y,1}\circ \Lambda_{Y,2}, 
&\Lambda_{X,2} &\;=\; [-1]_* (\Lambda_{Y,2}).
\end{align*}
Note that these are identities of cycles, including any possible multiplicities. Now take any irreducible subvariety $X_2\subset A$ such that $\Lambda_{X_2} \subseteq \Lambda_{X,2}$. By the above equations then 
\[ \Lambda_{Y_2} \;\subseteq\; \Lambda_{Y,2} \quad \textnormal{for} \quad Y_2 \;:=\; -X_2, 
\]
and looking at weights we find irreducible subvarieties $X_1, Y_1\subset A$ with $\Lambda_{X_1}\subseteq \Lambda_{X,1}$ and $\Lambda_{Y_1}\subseteq \Lambda_{Y,1}$ such that 
\begin{equation}
\label{eq:nZ}
 \Lambda_{nZ} \;\subseteq\; \Lambda_{X_1} \circ \Lambda_{Y_1} \;\subseteq\; \Lambda_{X_1} \circ \Lambda_{Y_1} \circ \Lambda_{X_2} \circ \Lambda_{Y_2}
\end{equation}
where the second inclusion comes from the antidiagonal $\Lambda_{\{0\}} \subseteq \Lambda_{X_2}\circ \Lambda_{Y_2}$. Now for $S=X,Y$ we have
\begin{equation}
\label{eq:S1S2}
 \Lambda_{S_1}\circ \Lambda_{S_2} \;\subseteq\; \Lambda_{S,1}\circ \Lambda_{S,2} \;=\; \cc(\delta_{nS}).
\end{equation}
In particular, any component of the left hand side maps to a subvariety of $A$ of dimension $\leq \dim(S)$. By~\eqref{eq:nZ} it follows that equality must hold for at least one such component, and then by looking at the right hand side of~\eqref{eq:S1S2} it follows that this component must be equal to $\Lambda_{nS}$. Since this component enters with multiplicity one, lemma~\ref{lem:multiplicity-one} implies
\[
 \dim(S_i) \;\leq\; \dim(S) \quad \textnormal{for} \quad i\;=\; 1,2.
\]
Taking $i=1$ and recalling that $\dim(X)+\dim(Y)=\dim(Z)<g$, we obtain from the first inclusion in~\eqref{eq:nZ} that
\[
 nZ \;=\; X_1 + Y_1
 \quad \textnormal{and} \quad 
 \dim X_1 \;=\; \dim X, \quad \dim Y_1 \;=\; \dim Y.
\]
In particular, for both $S=X,Y$ we have $\dim S_1 > 0$, hence $\deg(\Lambda_{S_1}) > 1$ and then a fortiori
\[
 \dim V_{S, 1} \;=\; \deg(\Lambda_{S,1}) \;\geq\; \deg(\Lambda_{S_1}) \;>\; 1.
\]
Since the $V_{S,1}$ are irreducible representations, it follows that they are nontrivial. \qed

\medskip 

The above argument should give more: It seems likely that $p: G_{XY} \twoheadrightarrow G_Z$ is an isogeny. For this we would like to argue that $\dim(X_2)=0$ by showing that one has an inclusion $\Lambda_{X_1+X_2} \subseteq \Lambda_{X_1}\circ \Lambda_{X_2}$. The problem is that a priori $\Lambda_{X_1 + X_2}$ could be negligible and in this case it is not taken into account by the clean product $\circ$ from example~\ref{ex:cleanconvolution}. However, in favorable cases this issue does not occur:

\begin{cor}
With notations as above, if the Gauss map $\gamma_{\, \bbP \Lambda}: \bbP \Lambda \to \bbP V$ for the projectivization of the characteristic variety $\Lambda = \Supp(\cc(\delta_Z))$  is a finite morphism, then the epimorphism $p: G_{XY} \twoheadrightarrow G_Z$ is an isogeny.
\end{cor}

{\em Proof.} This is shown in~\cite{KraemerJacobianThetaSummands} by computing Chern-Mather classes via lemma~\ref{lem:mather-ringhomo}; here the finiteness of the Gauss map means $\cc(\delta_Z) \in \scrL^{>g-1}(A)$.\qed

\medskip

\subsection{The adjoint representation} \label{sec:adjoint}

Any linear algebraic group has a distinguished representation, the adjoint representation on its Lie algebra. If $Z\subset A$ is a closed subvariety, let $\Ad_Z\in \Mhol(\scrD_A)$ be the unique clean module which corresponds to the adjoint representation of $G(\delta_Z)$. Combining the argument of~\cite[prop.~2.3]{KraemerMicrolocal} with the above we get \medskip

\begin{thm} \label{thm:summand} 
If $Z\subset A$ is irreducible with $\Stab(Z)=\{0\}$, then any
nontrivial decomposition into a sum of geometrically nondegenerate subvarieties
$Z = X + Y$ with $d_Z = d_X + d_Y$ satisfies \smallskip 
\[
  \min \{ d_X, d_Y\} \;\geq \; \delta \;=\; \tfrac{1}{2} \min 
  \bigl\{ \dim \Supp(\scrM) \mid \scrM \hookrightarrow \Ad_Z, \; \dim \Supp(\scrM) > 0 \bigr\}.  \medskip
\]
\end{thm} 

{\em Proof.} If $G$ is a connected semisimple group that is simple modulo its center, then for any nontrivial $U \in \Rep(G)$ the homomorphism $G\rightarrow \Gl(U)$ has finite kernel and so
\[
 Lie(G) \;\subseteq\; U\otimes U^* \;=\; \End(U).
\]
If $G$ is not simple modulo its center, the same argument still shows that $U\otimes U^*$ contains an irreducible summand of the adjoint representation. We apply this as follows: Up to an isogeny on the abelian variety we may assume all occuring Tannaka groups to be connected, so that the representations of their universal cover determine the corresponding clean $\scrD_A$-modules. For $S=X, Y$ the proposition says that
\[ \tilde{G}_Z \times \{1\}
 \quad \textnormal{acts nontrivially on } \quad  V_S.
\]
So by the above remarks on the adjoint representation there must be a nontrivial submodule
$\scrM \subseteq \Ad_Z$ with $\scrM\subseteq \delta_S*\delta_{-S}$.
Since for a connected semisimple group the adjoint representation has no one-dimensional summands, we must have $\dim \Supp(\scrM)>0$ and the claim follows because $\Supp(\scrM)\subseteq S-S$. \qed 

\medskip

If we fix a maximal torus, the nontrivial weights in the adjoint representation of a connected reductive group are by definition the roots of the group. The following criterion helps to estimate the support of $\Ad_Z$:

\begin{lem}
Let $Z\subset A$ be irreducible with $\Stab(Z)=\{0\}$. If the weights of~$V_Z$ with respect to a maximal torus contain some nonzero rational multiple of a root of~$G_Z$, then with notations as in theorem~\ref{thm:weyl-orbits}
\[ \dim \Supp(\Ad_Z) \;\geq\; \dim W \] 
for any  $\Lambda_W \subset h_*\cc(\delta_Z)$ in the image of the Weyl group orbit of this weight.
\end{lem}

{\em Proof.} The assumptions imply $[a]_*\Lambda_W \subseteq [b]_* \cc(\Ad_Z)$ for some $a,b\in \bbN$. \qed \bigskip


For ample divisors with at most rational singularities we know that any summand is geometrically nondegenerate~\cite[th.~1]{SchreiederDecomposable}, hence the above in particular applies to summands of theta divisors on indecomposable ppav's. The bound on the dimension of such summands is sharp whenever such summands are known to exist:
\begin{itemize} 
\item If $A$ is the Jacobian of a smooth projective curve $C$, then $\Ad_\Theta = \delta_{C-C}$. \smallskip 
\item If $A$ is the intermediate Jacobian of a smooth cubic threefold, 
$\Ad_\Theta = \delta_{\Theta}$. \smallskip
\end{itemize}
By the conjecture of Pareschi and Popa these should be the only theta divisors with summands of positive dimension. For curve summands this has been established in~\cite{SchreiederCurveSummands}. In higher dimension it remains open, but the above rules out various cases such as the following: 

\begin{thm} \label{thm:theta-not-a-sum}
Let $(A, \Theta)\in \scrA_g(\bbC)$. If $G(\delta_\Theta)$ is a symplectic group and~$\omega(\delta_\Theta)$ is its standard representation, then we have
$
 \dim \Supp(\Ad_\Theta) \geq g-1
$
and hence the theta divisor $\Theta \subset A$ cannot be written as a sum of positive-dimensional subvarieties.
\end{thm}

{\em Proof.} For symplectic groups the adjoint representation is the symmetric square of the standard representation. So twice any weight of the standard representation is a root, and the previous lemma shows $\dim \Supp(\Ad_\Theta) \geq g-1$. Thus any summand of the theta divisor must have dimension either zero or $\geq (g-1)/2$, which concludes the proof since $2\mid g$ if $\omega_u(\delta_\Theta)$ is symplectic. \qed

\medskip

\section{A criterion for almost simplicity} 

In order to control the adjoint modules from above and to compute the arising groups in more general cases, one needs a way to decide whether their Lie algebras are simple. We now give a sufficient criterion using characteristic cycles.

\subsection{Decomposable characteristic cycles} \label{sec:product}

We want to show that under the dictionary between weights and characteristic cycles in theorem~\ref{thm:weyl-orbits}, any nontrivial decomposition of the Lie algebra of the Tannakian group induces a decomposition on characteristic cycles. Let $\deg: \scrL(A) \to \bbZ$ be the homomorphism sending a clean subvariety $\Lambda \subset T^*A$ to the degree $\deg(\Lambda)=\deg(\gamma_\Lambda)$ of its Gauss map.

\begin{prop} \label{prop:product}
Let $G=G(\scrM)$ for a simple holonomic module $\scrM \in \Mhol(\scrD_A)$ 
with $\Stab(\scrM)=\{0\}$. If the Lie algebra $\mathrm{Lie}(G)$ 
is not simple modulo its center, then
\[
 [n]_*
 \cc(\scrM) \;=\; \Lambda_{1} \circ \Lambda_{2}
\]
for some $n\in \bbN$ and effective clean cycles $\Lambda_i\in \langle \cc(\scrM)\rangle$ of degree $\deg(\Lambda_i)>1$.
\end{prop}

{\em Proof.} Under the assumptions of the proposition, the connected component~$G^\circ$ is not simple modulo its center $Z=Z(G^\circ)$, and by the theory of connected reductive groups we can find nontrivial connected semisimple groups $G_1, G_2\neq \{1\}$ and an epimorphism
\[
 \widehat{G} \;=\; G_1\times G_2 \times Z \; \twoheadrightarrow \; G^\circ
\]
whose kernel is a finite central subgroup. As $\Stab(\scrM)=\{0\}$, we know that~$\omega_u(\scrM)$ restricts to an irreducible representation of the connected component~\cite[cor.~1.4]{KraemerMicrolocal}, so
$
 U = \omega_u(\scrM)|_{\tilde{G}} 
$
is an irreducible representation of the covering group. As such it splits as a product
$
 U \simeq U_1\boxtimes U_2$ of irreducible
$ U_1 \in \Rep(G_1),  U_2 \in \Rep(G_2\times Z)$,
and the claim then follows from theorem~\ref{thm:inverse-galois} and remark~\ref{rem:effective}. \qed

\subsection{A criterion for divisors}

Now let $Z\subset A$ be a divisor and $p: \tilde{Z} \to Z$ a resolution of singularities. The universal property of the Albanese variety gives a unique $\hat{p}$ making the following diagram commute:
\[
\xymatrix@M=0.5em{
 \tilde{Z} \ar[r]^-{alb} \ar[d]_-p & Alb(\tilde{Z}) \ar@{..>}[d]^-{\exists ! \, \hat{p}} \\
 Z \ar@{^{(}->}[r] & A 
}
\]
We want $\hat{p}$ to be an isomorphism. As the Albanese variety is a birational invariant, this condition does not depend on the chosen resolution and we abbreviate it by writing $A=Alb(Z)$. The following criterion applies in particular to $\scrM= \delta_Z$ in many cases:

\begin{thm} \label{thm:simple}
Let $Z\subset A$ be a geometrically nondegenerate irreducible reduced symmetric divisor such that $\Stab(Z)=\{0\}$ and $Alb(Z)=A$. Let $\scrM\in \Mhol(\scrD_A)$ be simple and
\[
 \cc(\scrM) \; = \; \Lambda_Z + \sum_{d_W < d_Z} c_W\cdot \Lambda_W
 \quad \textnormal{\em with} \quad c_W \in \bbN_0.
\]
Then $\mathrm{Lie}(G(\scrM))$ is simple modulo its center in each of the following situations: \medskip
\begin{enumerate}
\item If $\deg(\Lambda_Z) > \tfrac{1}{3} \deg(\cc(\scrM))$. \smallskip
\item If $n_W\neq 0$ at most for $\dim(W)=0$. \smallskip
\item If $\gamma_Z: \bbP\Lambda_Z \to \bbP V$ is finite and $[2m]_*\cc(\scrM) \neq \Lambda_{mZ} \circ \Lambda_{mZ}$ for all $m\in \bbN$. \smallskip
\item If $[n]_*\cc(\scrM)$ is reduced for all $n\in \bbN$. \medskip
\end{enumerate}
\end{thm}

{\em Proof.} If the connected component is not simple modulo its center, pick $n\in \bbN$ and $\Lambda_i \in \scrL(A)$ with $\deg(\Lambda_i) > 1$ as in proposition~\ref{prop:product}. Since $\Stab(Z)=\{0\}$, the map 
\[
 [n]: \quad Z \;\longrightarrow\; A, \quad z\mapsto nz 
\]
is birational onto its image. So our assumption implies that the conormal variety to this image $Y=nZ$ enters the clean conic Lagrangian cycle $[n]_*\cc(\scrM) = \Lambda_1 \circ \Lambda_2$ with multiplicity one, and all other components of this cycle are conormal varieties to subvarieties of codimension $>1$ in $A$. Let $\Lambda_{Z_i} \subset \Supp(\Lambda_i)$ be the unique irreducible components with 
\[ \Lambda_Y \;\subseteq \;\Supp(\Lambda_{Z_1} \circ \Lambda_{Z_2}), \]
and write $\Lambda_i = \Lambda_{Z_i} + \Lambda_i'$ with effective cycles $\Lambda_i'\in \scrL(A)$, possibly zero. Then we have
\begin{eqnarray}
 \Lambda_{Z_1} \circ \Lambda_{Z_2}
 & \;=\; & \Lambda_Y + \cdots \label{eq:smallconv} \\
 (\Lambda_{Z_1} + \Lambda_1') \circ (\Lambda_{Z_2} + \Lambda_2')
 & \;=\; &
 \Lambda_Y + \cdots  \label{eq:bigconv} \\
 \deg(\Lambda_{Z_1}+\Lambda_1') \cdot \deg(\Lambda_{Z_2}+\Lambda_2') 
 &\;=\; & \deg(\cc(\scrM)) \label{eq:degrees}
\end{eqnarray}
where $\cdots$ stands for a nonnegative linear combination of conormal varieties to subvarieties of codimension $>1$ in $A$.
We now argue by contradiction. \medskip 

First we claim that $\dim(Z_i) >0$ for both $i = 1,2$. If not, then by symmetry we can assume $Z_2=\{p\}$ is reduced to a single point. Then $Z_1 = Y-p$ by~\eqref{eq:smallconv}, and by~\eqref{eq:bigconv} it follows that $\Supp(\Lambda_2')$ cannot contain the conormal variety to any point on the abelian variety. Since conormal varieties to points are the only ones whose Gauss map has degree one, we get $\deg(\Lambda_2')\geq 2$. But $\deg(\Lambda_{Y-p})=\deg(\Lambda_Z)$, so case (1) is excluded by~\eqref{eq:degrees}.
In case (2) any irreducible component of the convolution $\Lambda_{Y-p}\circ \Lambda_2' \neq 0$ would be supported over a point. By the last part of lemma~\ref{lem:mindim-of-convolution} and the symmetry of $Y\subset A$ then $\Lambda_{Y+q}\subseteq \Lambda_2'$ for some $q\in A(\bbC)$, but then the diagonal in the fibered product would give an irreducible component of the form
$
 \Lambda_{2Y+q-p} \subseteq \Supp(\Lambda_{Y-p}\circ \Lambda_2')
$
contradicting~\eqref{eq:smallconv}-\eqref{eq:bigconv}. 
Case (3) would result in 
\[
 [Y] \;=\; c_{M,g-1}(\Lambda_1\circ \Lambda_2) \;\geq \; \deg(\Lambda_2) \cdot [Y] \;>\; [Y]
\]
by~\eqref{eq:bigconv} and lemma~\ref{lem:mather-ringhomo}. In the remaining case (4) the right hand side of~\eqref{eq:bigconv} is a reduced cycle, hence so is the convolution product $\Lambda_{Z_1} \circ \Lambda_2' \neq 0$. Since $Z_1\subset A$ is a divisor, lemma~\ref{lem:multiplicity-one} then forces $\Lambda_{Z_1} \circ \Lambda_2'$ to be a sum of conormal varieties to divisors on $A$. But this is not the case for any of the summands $\cdots$  in~\eqref{eq:bigconv}.

\medskip

So $\dim(Z_i)\geq 1$ for both $i = 1,2$. By~\eqref{eq:smallconv} and lemma~\ref{lem:multiplicity-one} we have dominant rational maps
\[
 f_i: \quad 
 Z \;\stackrel{[n]}{\twoheadrightarrow} \; Y \;\dashrightarrow\; Z_i
\]
for $i=1,2$.
By the extension property of rational maps from smooth projective varieties to abelian varieties~\cite[th.~3.1]{MilneAV}, these maps can be extended to surjective morphisms 
\[
  \tilde{f}_i: \quad
  \tilde{Z} \; \twoheadrightarrow\; Z_i \;\subset \; A
\]
where $\tilde{Z} \twoheadrightarrow Z$ is a resolution of singularities. Our assumption $A=Alb(Z)$ then implies that the original rational maps $f_i$ are defined everywhere and fit into a commutative diagram
\[
\xymatrix@M=0.6em@C=2em@R=2em{
 \tilde{Z} \ar[r] \ar[d]_-{alb} & Z  \ar[r]^{f_i} \ar@{^{(}->}[d] & Z_i \ar@{^{(}->}[d] \\
 Alb(Z) \ar@{=}[r] & A \ar@{..>}[r]^{\exists ! g_i} & A
}
\]
where $g_i$ comes from the universal property of the Albanese variety. Replacing~$Z_i$ by a translate we may assume $Z_i$ is symmetric and $g_i$ is a homomorphism. Then the image $B_i=g_i(A) \subseteq A$ is an abelian subvariety. On the other hand $Z_i \subset A$ cannot be a translate of an abelian subvariety, indeed $\deg(\Lambda_{Z_i})>0$. So the geometric nondegeneracy of $Z$ and the above diagram imply that the map $f_i: Z\twoheadrightarrow Z_i$ is generically finite, and 
$
 g_i:  A \twoheadrightarrow B_i = A
$
must then be an isogeny for dimension reasons. Since $\Stab(Z)=\{0\}$, it follows that $f_i: Z \rightarrow Z_i$ is birational. 

\medskip

Note that since we have adjusted our translates such that the $g_i: A\to A$ are homomorphisms, the symmetry of $Z$ implies that $Z_i$ is symmetric. On the other hand
$f_i \neq -f_i$
since $Z_i \not\subset A[2]$. It follows that in the fiber product $\bbP \Lambda_{Z_1} \times_{\bbP V} \bbP \Lambda_{Z_2}$ we have components 
\[
 \Delta^\pm \;=\; \{ (f_1(z), \pm f_2(z), \xi) \in A\times A \times \bbP V \mid (z, \xi) \in \Lambda_Z \} 
\]
and these two components are distinct since the $f_i$ are birational. Both components are dominant over $\bbP V$ and so their image under the addition morphism $\varpi$ gives two components
\[
 \Lambda_{Y^\pm} \;\subset\; \Lambda_{Z_1} \circ \Lambda_{Z_2}
 \quad \textnormal{where} \quad
 Y^\pm \;=\; (g_1\pm g_2)(Z).
\]
Here $Y^+ = Y$ because $g_1+g_2=[n]$. Looking at the multiplicities of the components in equation~\eqref{eq:smallconv} we then obtain $\dim Y^- < \dim Z$. So $g_1-g_2: Z\twoheadrightarrow Y^-$ is not generically finite, and since $Z$ is geometrically nondegenerate, it follows that the image
$
 Y^- =(g_1-g_2)(A)
$
is an abelian subvariety. But $\deg(\Lambda_{Y^-}) > 0$, so $Y^- = \{0\}$ must be a single point. Thus $g_1 = g_2$. Since $f_1+f_2 = [n]$, it follows by looking at the induced map on the universal cover of the abelian variety that $n=2m$ is even and $g_1=g_2=[m]$, so
\[ 
  Z_1 \;=\; Z_2 \;=\; mZ.
\]
To show that this is impossible in each of the four cases listed in the theorem, recall that
\[
 [n]_*\cc(\scrM) \;=\; \Lambda_1\circ \Lambda_2
 \;=\; (\Lambda_{Z_1} + \Lambda_1') \circ (\Lambda_{Z_2} + \Lambda_2'). \medskip
\]
So case (4) is excluded by observing that this cycle contains the component $\Lambda_{\{0\}}$ with multiplicity at least $\deg(\varpi: \Delta^- \twoheadrightarrow \Lambda_{\{0\}}) = \deg(\Lambda_Z)\geq 2$. The remaining three cases can be excluded by a computation of Chern-Mather classes: In degree zero we get
$
 \deg(\Lambda_Z)^2 \leq \deg(\cc(\scrM)),
$
so case (1) could only happen if $\deg(\cc(\scrM)) \leq 3$ and this is irrelevant for us since any reductive Lie algebra with a faithful representation of dimension at most three is simple modulo its center. In degree one lemma~\ref{lem:mather-ringhomo} implies 
\[
2m^2 \deg(\Lambda_{Z})c_{M,1}(\Lambda_{Z}) \;\leq\; n^2 \, c_{M,1}(\cc(\scrM)),
\]
where the inequality means that the right hand is equal to the left hand side plus an effective cycle. In case (2) we would have $c_{M,1}(\cc(\scrM))=c_{M,1}(\Lambda_Z)$ and it would follow from the above inequality that $\deg(\Lambda_Z)=2$. Then the Gauss map for the symmetric divisor $Z\subset A$ would induce a birational map $Z/\langle \pm \id_A\rangle \dashrightarrow \bbP V$ which is impossible for the same reason as in~\cite[lemma~3.1]{CGS}: Ran's definition of geometric nondegeneracy in~\cite[p.~466]{RanSubvarieties} implies that on any geometrically nondegenerate divisor~$Z\subset A$ there are non-zero holomorphic $2$-forms invariant under $-\id_A$. So it only remains to exclude the case (3). Taking Chern-Mather classes in degree $g-1$ in~\eqref{eq:bigconv} we get
\[
c_{M,g-1}\bigl(\Lambda_{mZ}\circ \Lambda_2' + \Lambda_1'\circ \Lambda_{mZ} + \Lambda_1'\circ \Lambda_2' \bigr) \;=\; 0.
\]
In case (3) the finiteness of the Gauss map allows to use lemma~\ref{lem:mather-ringhomo}. Effective clean cycles have effective Chern-Mather classes and the Pontryagin product preserves effectivity, so we get
$
\deg(\Lambda_1') \cdot c_{M,g-1}(\Lambda_{mZ})  \;=\; \deg(\Lambda_2') \cdot c_{M, g-1}(\Lambda_{mZ}) \;=\; 0.
$
This implies $\Lambda_1'=\Lambda_2'=0$ because we are dealing with clean effective cycles. Altogether it would follow that
$[n]_* \cc(\scrM) \;=\; \Lambda_{mZ}\circ \Lambda_{mZ}$
which we had excluded in (3). \qed

\begin{cor} \label{cor:minuscule}
For $\scrM$ as in theorem~\ref{thm:simple} with $[n]_* \cc(\scrM)$ reduced for all $n\in \bbN$,  we have 
\[ G \;\subseteq \; G(\scrM)^\circ \;\subseteq \; \bbG_m \cdot G \;\subseteq\; \Gl(W), \smallskip \]
where $(G, W)$ is one of the pairs that are listed in tables~\ref{tab:minuscule},~\ref{tab:wmf} of the appendix.
\end{cor}

{\em Proof.} By theorem~\ref{thm:simple} the connected component $G(\scrM)^\circ$ is simple modulo its center, and since $\scrM \in \Mhol(\scrD_A)$ is simple with trivial stabilizer $\Stab(\scrM)=\{0\}$, this connected component acts irreducibly on $W=\omega(\scrM)$. So by Schur's lemma its center acts via a scalar, which means that we have inclusions
$G \subseteq  G(\scrM)^\circ  \subseteq \bbG_m \cdot G$ where $G$ is the derived group of the connected component. Our assumption on the reducedness of $[n]_*\cc(\scrM)$ for all $n\in \bbN$ implies by theorem~\ref{thm:weyl-orbits} that $W$ is weight multiplicity free, so the claim follows from lemma~\ref{lem:wmf}. \qed

\subsection{Application to theta divisors}

By~\cite{EL} theta divisors on indecomposable ppav's are normal and irreducible so that the Albanese assumption from above is satisfied:

\begin{lem} \label{lem:ppav}
For any indecomposable ppav $(A, \Theta) \in \scrA_g(\bbC)$ of dimension $g>1$, the theta divisor satisfies $A=Alb(\Theta)$.
\end{lem}

{\em Proof.} Let $p: \tilde{\Theta} \to \Theta$ be a resolution of singularities, and put $B=Alb(\tilde{\Theta})$. We have a diagram
\[
\xymatrix@M=0.5em{
 \tilde{\Theta} \ar[r]^-{alb} \ar[d]_-p & B \ar@{..>}[d]^-{\exists ! \, \hat{p}} \\
 \Theta \ar@{^{(}->}[r] & A
}
\]
where $\hat{p}$ is an isogeny~\cite[rem.~3.4]{EL}. As~$p$ is birational, we know
$
 D \;=\; alb(\tilde{\Theta}) \;\subset\; B
$
is an irreducible divisor which is not stable under any translation: If $e\in \Stab(D)$, then $\hat{p}(e) \in \Stab(\Theta)=\{0\}$, so 
$\hat{p}(d) = \hat{p}(d+e)$
for all $d\in D(\bbC)$; since $d, d+e\in D$, the birationality of
\[
 p: \quad 
 \xymatrix@M=0.5em{
 \tilde{\Theta} \ar@{->>}[r]^{alb} & D \ar@{->>}[r]^-{\hat{p}} & \Theta
 }
\]
then forces $e=0$. Hence $\Stab(D)=\{0\}$ as claimed. Writing the preimage of $\Theta$ as a union
\[
 \hat{p}^{-1}(\Theta) \;\;=\; \bigcup_{e\in \ker(\hat{p})} (D + e),
\]
we get that on the right hand side the irreducible divisors $D+e \subset B$ are pairwise distinct. All these divisors are ample since they are all numerically equivalent and their union $\hat{p}^{-1}(\Theta) \subset B$ is an ample divisor. Hence if there exists $e\in \ker(\hat{p}) \setminus \{0\}$, it follows that
$ \varnothing \neq D\cap (D+e) \subset D $
is nonempty of codimension one, in which case the preimage $\hat{p}^{-1}(\Theta)$ will be singular in codimension one. But this is impossible since $\hat{p}$ is an isogeny so that the preimage of any normal divisor must be normal.
\qed 

\medskip

We can thus apply theorem~\ref{thm:simple} to the module $\scrM=\delta_\Theta\in \Mhol(\scrD_A)$. Notice that while the polarization determines the theta divisor only up to a translate, the connected semisimple group
\[
 G_\Theta \;=\; [G(\delta_\Theta)^\circ, G(\delta_\Theta)^\circ]
\]
is independent of the chosen translate due to lemma~\ref{lem:translate}. Thus we get:

\begin{thm} \label{thm:theta-multiplicity-free}
Let $(A, \Theta) \in \scrA_g(\bbC)$ be a ppav 
such that $[n]_*\cc(\delta_\Theta)$ is reduced for all $n\in \bbN$. Then the connected semisimple group $G_\Theta$ and its irreducible faithful representation $W=\omega(\delta_\Theta)$ appear in table~\ref{tab:minuscule} or~\ref{tab:wmf}.
\end{thm}

{\em Proof.} Up to a translation we may assume that the theta divisor is symmetric, so that the representation $\omega(\delta_\Theta)$ is isomorphic to its dual. Since the restriction of this representation to the connected component $G(\delta_\Theta)^\circ$ remains irreducible, it follows from Schur's lemma that this connected component must be semisimple and hence equal to $G_\Theta$. So the claim follows from corollary~\ref{cor:minuscule}. \qed

\medskip

This contains theorem~\ref{thm:theta_odp} from the introduction: If $\Theta$ is smooth except for~$k$ ordinary double points, then the classical Gauss map has degree
$ \deg(\Lambda_\Theta) =  g! - 2k $
by~\cite{CGS}; on the other hand $\dim(\omega(\delta_\Theta)) = \deg(\cc(\delta_\Theta))$ is always the degree of the characteristic cycle, and the difference between the two is given in lemma~\ref{lem:cc}.

\subsection{Stratifications of moduli spaces} \label{sec:stratification}

It follows from~\cite[prop.~7.4]{KrWSchottky} that the assignment
$
 (A, \Theta) \mapsto (G(\delta_\Theta), \omega(\delta_\Theta))
$
defines a constructible stratification of the moduli space of ppav's. This stratification is related to the stratification by the degree of the Gauss map in~\cite{CGS} but includes finer information. Let $\scrA_g^{ind}\subseteq \scrA_g$ be the locus of indecomposable ppav's. For $g=4$ we have
\[
 \scrA_4^{ind} \;=\;
 \scrA_4^{sm}
 \;\sqcup\;
 \scrJ_4^{nh}
 \;\sqcup\;
 \scrJ_4^{h}
 \;\sqcup\; \bigsqcup_{k=1}^{10} \,
 \Theta_{null,4}^{k}
\]
where \medskip
\begin{itemize}
\item $\scrA_4^{sm}$ \hspace*{0.4em} is the locus of ppav's with a smooth theta divisor, \smallskip 
\item $\scrJ_4^{nh}$ \hspace*{0.3em} is the locus of Jacobians of nonhyperelliptic curves, \smallskip 
\item $\scrJ_4^h$ \hspace*{0.8em} is the locus of Jacobians of hyperelliptic curves, \smallskip 
\item $\Theta_{null, 4}^k$ is the locus of non-Jacobians with precisely $k$ vanishing thetanulls. \medskip
\end{itemize}
Indeed, the only singularities of indecomposable theta divisors on non-Jacobian abelian fourfolds are ordinary double points given by vanishing thetanulls~\cite{SmithVarleyComponents}; the maximum number of $k=10$ such vanishing theta nulls is obtained for a unique abelian fourfold discovered by Varley~\cite{VarleyWeddlesSurfaces}~\cite{DebarreAnnulation}. Table~\ref{tab:ppavs} lists  \medskip
\begin{enumerate} 
\item the degree $\deg(\Lambda_\Theta)$ of the classical Gauss map, \smallskip
\item the degree $\deg(\cc(\delta_\Theta)) = \dim(\omega(\delta_\Theta))$ of the characteristic cycle, \smallskip
\item the representation $\omega(\delta_\Theta)$ of the reductive group $G(\delta_\Theta)$, \medskip
\end{enumerate} 
where for the latter we denote by $\varpi_1, \varpi_2, \dots$ the fundamental representations.

\begin{table}[h] 
\[
\small
\def\arraystretch{1.4}
\begin{array}{|l||l|l|l|l|} \hline 
  & \deg(\gamma_\Theta) & \dim(\omega(\delta_\Theta)) & \omega(\delta_\Theta) & G(\delta_\Theta)  \\ \hline 
 \scrA_4^{sm} & 24 & 24 & \varpi_1 & \Sp_{24}(\bbC) \\
 \scrJ_4^{nh} & 20 & 20 & \varpi_3 & \Sl_6(\bbC)/\mu_3  \\
 \scrJ_4^h  & 8 & 14 & \varpi_3 & \Sp_6(\bbC)  \\
 \Theta_{null, 4}^k & 24 - 2k & 24 - 2k  & \varpi_1 & \Sp_{24-2k}(\bbC) \; \textnormal{(assuming $\gamma_\Theta$ is finite if $k=4$)} \\  \hline 
\end{array} 
\]
\caption{Invariants of principally polarized abelian fourfolds}
\label{tab:ppavs}
\end{table}

Note that for $k=2$ neither the degree of the classical Gauss map nor the one of the characteristic cycle suffices to characterize Jacobians among all ppav's, but the pair $(G(\delta_\Theta),\omega(\delta_\Theta))$ does. Let us briefly check the new part of the table:

\begin{lem} 
For $k\neq 2$, theorem~\ref{thm:theta_odp} implies that every ppav $(A, \Theta)\in \Theta_{null, 4}^k(\bbC)$ has
\[
G(\delta_\Theta) \simeq \Sp_{24-2k}(\bbC) 
\quad \textnormal{\em and} \quad 
\omega(\delta_\Theta) \simeq \varpi_1. \medskip
\]
For $k= 2$ this remains true at least if the Gauss map $\gamma_\Theta: \bbP \Lambda_\Theta \to \bbP V$ is finite. 
\end{lem}

{\em Proof.} For $k\neq 2$ this is clear by direct inspection. For $k = 2$ the only possible alternative would be that $G(\delta_\Theta)\simeq \Sl_6(\bbC)/\mu_3$, $\omega(\delta_\Theta) \simeq \varpi_3$. In this case theorem~\ref{thm:fake-jacobian} shows 
\[
 [3]_*\cc(\delta_\Theta) \;=\; \Alt^3(\Lambda)
 \quad \textnormal{where $\Lambda \in \scrL(A)$ is effective up to isogeny}.
\]
If the Gauss map is finite, a computation of Chern-Mather classes via lemma~\ref{lem:mather-ringhomo} gives 
\[
 c_{M,i}(\Lambda) \;=\; 
 \begin{cases}
 6 & \textnormal{for $i=0$}, \\
 \tfrac{2}{3} [\Theta]^3 & \textnormal{for $i=1$}, \\
 0 & \textnormal{for $i>1$},
 \end{cases}
\]
by example~\ref{ex:elementary-symmetric} and corollary~\ref{cor:chern-gauss}. For any $n\in \bbN$ such that the cycle $[n]_*(\Lambda)$ is effective, it then follows from lemma~\ref{lem:mather-transversal} that this cycle must be supported over a curve $Z\subset A$. Then the image of the theta divisor under the isogeny $[3n]$ is a sum of copies of this curve. Taking preimages under this isogeny we get that the theta divisor is a sum of curves, so $(A, \Theta)$ is a Jacobian by~\cite{SchreiederCurveSummands}. \qed

\medskip

Finally, let us give an example to illustrate how the dictionary between weights and characteristic cycles can be used beyond theorem~\ref{thm:theta-multiplicity-free}. We only give the simplest case but the method obviously generalizes:

\begin{lem} \label{lem:torsion-difference}
If $(A, \Theta) \in \scrA_5(\bbC)$ is a ppav with a symmetric theta divisor that is smooth except for two distinct ordinary double points $e_1, e_2$, then one of the following two cases occurs: \smallskip 

\begin{enumerate} 
\item Either $e_1,e_2\in A[2]$, in which case $G(\delta_\Theta) \simeq O_{118}(\bbC)$. \smallskip 
\item Or $e_1=-e_2\notin A[2]$, in which case $G(\delta_\Theta) \simeq \SO_{118}(\bbC)$. \smallskip
\end{enumerate}
\end{lem}

{\em Proof.}
Since the theta divisor is symmetric and the dimension $g=5$ is odd, the representation $\omega(\delta_\Theta)$ is orthogonal. By lemma~\ref{lem:cc} the characteristic cycle is given by
\[ \cc(\delta_\Theta) \;=\; \Lambda_\Theta + \Lambda_{e_1} + \Lambda_{e_2}. \] 
If the weights of the representation form a single Weyl group orbit, then we are done by lemma~\ref{lem:wmf} because the distinction between the orthogonal and special orthogonal group can be read off from $\det(\delta_\Theta)=\delta_{e_1+e_2}$. So it only remains to exclude the case that there are more than one Weyl group orbits. Since one of them must be of size $\deg(\Lambda_\Theta)$, the remaining orbits would have to be of size at most two. But for any simple complex Lie algebra, any nontrivial Weyl group orbit has size at least the rank of the Lie algebra, which in our case must be bigger than two as otherwise the first orbit could not be so big. So all nontrivial weights are in the same Weyl group orbit, i.e.~the representation is quasi-minuscule. But these are all known and there is no such of dimension $\deg(\cc(\delta_\Theta))=118$. \qed

\medskip

\subsection{Appendix: Multiplicity free representations}

\label{sec:multiplicity-free}


Although the dictionary between Weyl group orbits and characteristic cycles works in general, the simplest situation is when $\scrM$ is {\em essentially multiplicity free} in the sense that the clean cycle~$[n]_*(\cc(\scrM))$ is reduced for all $n\in \bbN$, see~\cite{KraemerMicrolocal}. Then any Weyl group orbit of weights that enters the corresponding representation must do so with multiplicity one, i.e.~$\omega_u(\scrM)$ is a ~{\em weight multiplicity free} representation. For the simple Dynkin types there are only very few such representations, most of them are minuscule in the sense that their weights form a single orbit under the Weyl group:

\begin{lem} \label{lem:wmf}
If $G$ is a connected semisimple group which is simple modulo its center and acts faithfully on some weight multiplicity free irreducible $W\in \Rep(G)$, then $(G, W)$ appears on table~\ref{tab:minuscule} or~\ref{tab:wmf} below.  
\end{lem}

{\em Proof.}
The weight multiplicity free irreducible representations of the simple Lie algebras are classified in~\cite[th.~4.6.3]{HowePerspectives}~\cite[sect.~A]{StembridgeMultiplicity}, and the tables list the images of the corresponding simply connected groups under these representations. Whether an irreducible representation is orthogonal, symplectic or not self-dual can be read off from its highest weight~\cite[sect.~3.2.4 and exercise~9]{GW}.
\qed

\begin{table}[h] 
	\[
	\small
	\def\arraystretch{1.4}
	\begin{array}{|r||r|r|r|r|r|r|r|} \hline
	G & \Sl_{n}/\mu_{(k, n)} & Spin_{2n+1} & \Sp_{2n} & \SO_{2n} & Spin_{2n} & E_6 & E_7 \\ \hline 
	\dim W & {n \choose k} \; \textnormal{for} \; 1\leq k\leq n & 2^n & 2n  & 2n & 2^{n-1} & 27 & 56 \\ \hline 
	\textnormal{symplectic?} & n = 2k \notin 4\bbZ  & n\equiv 1,2 \, (4) & \textnormal{yes} & \textnormal{no} & n\equiv 2 \, (4) & \textnormal{no} & \textnormal{yes} \\
	\textnormal{orthogonal?} & n = 2k \in 4\bbZ & n\equiv 0,3 \, (4) & \textnormal{no} & \textnormal{yes} & n\equiv 0 \, (4) & \textnormal{no} & \textnormal{no} \\ \hline
	\end{array}
	\medskip
	\]
	\caption{The minuscule representations $W$ and their images $G \subset \Gl(W)$ for the simple Dynkin types. All these are fundamental representations. In the last two rows we list whether they are symplectic, orthogonal or not self-dual.}
	\label{tab:minuscule}
\end{table}
\vspace*{-2.5em}
\begin{table}[h]
	\[
	\small
	\def\arraystretch{1.4}
	\begin{array}{|r||r|r|r|r|r|r|r|} \hline
	G & \Sl_{n}/\mu_{(k, n)} & \SO_{2n+1} & \Sp_6 & G_2 \\ \hline 
	\dim W &  {n+k-1 \choose k} \; \textnormal{for} \; k > 1 & 2n+1 & 14 & 7 \\ \hline 
	\textnormal{symplectic?} & \textnormal{no} & \textnormal{no} & \textnormal{yes} & \textnormal{no} \\
	\textnormal{orthogonal?} & \textnormal{no} & \textnormal{yes} & \textnormal{no} & \textnormal{yes} \\ \hline
	\end{array}
	\medskip
	\]
	\caption{The weight multiplicity free nonminuscule representations $W$ and their images $G\subset \Gl(W)$ for the simple Dynkin types. All these are fundamental representations except for the symmetric powers $W=Sym^k(\bbC^n)\in \Rep(\Sl_n)$.}
	\label{tab:wmf}
\end{table}

\medskip 

\begin{lem} \label{lem:cc}
If $Z\subset A$ is a reduced divisor whose singular locus consists of finitely many ordinary double points $e_i$, then \smallskip
\[
 \CC(\delta_Z) \;=\;
 \begin{cases}
  \Lambda_Z & \textnormal{\em if $g$ is even}, \\
 \Lambda_Z + \sum_{i} \Lambda_{e_i} & \textnormal{\em if $g$ is odd}.
 \end{cases}
\]
So $\delta_Z$ is essentially multiplicity free iff \smallskip 
\begin{itemize}
\item $\Stab(Z)=\{0\}$, and \smallskip
\item no two of the $e_i$ differ by a torsion point if $2\nmid g$. \smallskip
\end{itemize}
\end{lem}

{\em Proof.} The characteristic cycle can be computed locally in the classical topology, so we can apply the result of~\cite[th.~4.2]{NTCharacteristic}: If $X\subseteq \bbC^g$ is an open neighborhood of the origin and $f: X\to \bbC$ is a holomorphic function such that $Z = f^{-1}(0)$ is smooth outside the origin, then
$
 \CC(\delta_Z) = \Lambda_Z + (\mu_0 - N_1) \cdot \Lambda_0
$
where $\mu_0$ is the Milnor number of a generic hyperplane section of $Z$ through the origin and where $N_1$ denotes the number of Jordan blocks for the eigenvalue one in the Milnor monodromy on $H^{g-1}(f^{-1}(t), \bbC)$ for small $t\neq 0$. For ordinary double points $\mu_0 = 1$, and the local Picard-Lefschetz formula says that the Milnor monodromy acts on the one-dimensional space $H^{g-1}(f^{-1}(t), \bbC)$ by $(-1)^g$.
\qed

\bigskip \bigskip

{\em Acknowledgements.} I would like to thank Giulio Codogni for valuable comments on a previous version of the paper, and Michael Dettweiler, Gavril Farkas, Javier Fres\'an, Bert van Geemen, Bruno Klingler, Mihnea Popa and Stefan Schreieder for stimulating discussions on various related topics.

\bigskip \medskip

\bibliographystyle{amsplain}
\bibliography{Bibliography}

\end{document}